\documentclass[11pt,letterpaper]{amsart}
\usepackage[latin1]{inputenc}
\usepackage[francais,english]{babel}    
\usepackage{amssymb}
\usepackage{amsmath}
\usepackage{amsthm}
\usepackage{multicol}
\usepackage{graphicx}
\usepackage{color}
\usepackage[T1]{fontenc}
\usepackage{yfonts}
\usepackage{placeins}
\usepackage{array}
\usepackage{dsfont}
\usepackage{stmaryrd}
\usepackage{verbatim}
\usepackage{caption}
\usepackage{transparent}
\usepackage{bbold}
\usepackage{hyperref}
\usepackage{enumitem}
\usepackage{tikz}
\usepackage{float}
\usetikzlibrary{patterns}
\usetikzlibrary{arrows.meta}
\usetikzlibrary{decorations.markings}
\usetikzlibrary{matrix,arrows,calc,fit,cd,positioning,intersections,arrows.meta,braids}
\usepackage{pgfplots}
\usepackage{tikz-3dplot}

\newtheorem{thm}{Theorem}[section]
\newcommand{\bthm}{\begin{thm}}
\newcommand{\ethm}{\end{thm}}

\newtheorem{thmi}{Theorem}
\newcommand{\bthmi}{\begin{thmi}}
\newcommand{\ethmi}{\end{thmi}}

\newtheorem{cori}[thmi]{Corollary}
\newcommand{\bcori}{\begin{cori}}
\newcommand{\ecori}{\end{cori}}

\newtheorem{mthm}{Theorem}
\newcommand{\bmthm}{\begin{mthm}}
\newcommand{\emthm}{\end{mthm}}

\newtheorem{mcor}[mthm]{Corollary}
\newcommand{\bmcor}{\begin{mcor}}
\newcommand{\emcor}{\end{mcor}}
\newtheorem{mconj}[mthm]{Conjecture}
\newcommand{\bmconj}{\begin{mconj}}
\newcommand{\emconj}{\end{mconj}}
\newtheorem{mpro}[mthm]{Proposition}
\newcommand{\bmpro}{\begin{mpro}}
\newcommand{\empro}{\end{mpro}}

\newtheorem*{conj}{Conjecture}
\newcommand{\bconj}{\begin{conj}}
\newcommand{\econj}{\end{conj}}

\newtheorem*{question}{Question}
\newcommand{\bq}{\begin{question}}
\newcommand{\eq}{\end{question}}

\newtheorem*{thn}{Theorem}
\newcommand{\bthn}{\begin{thn}}
\newcommand{\ethn}{\end{thn}}

\newtheorem{exo}{Exercise}
\newcommand{\bex}{\begin{exo}}
\newcommand{\eex}{\end{exo}}

\newtheorem{sol}{Solution}
\newcommand{\bsol}{\begin{sol}}
\newcommand{\esol}{\end{sol}}

\newtheorem{pro}[thm]{Proposition}
\newcommand{\bpro}{\begin{pro}}
\newcommand{\epro}{\end{pro}}

\newtheorem{cor}[thm]{Corollary}
\newcommand{\bcor}{\begin{cor}}
\newcommand{\ecor}{\end{cor}}

\newtheorem{lem}[thm]{Lemma}
\newcommand{\blem}{\begin{lem}}
\newcommand{\elem}{\end{lem}}

\theoremstyle{definition}

\newtheorem{defi}[thm]{Definition}
\newcommand{\bdf}{\begin{defi}}
\newcommand{\edf}{\end{defi}}

\newtheorem*{defis}{Definition}
\newcommand{\bdfs}{\begin{defis}}
\newcommand{\edfs}{\end{defis}}

\newtheorem*{rmk}{Remark}
\newcommand{\brk}{\begin{rmk} \upshape}
\newcommand{\erk}{\end{rmk}}

\newtheorem*{rmks}{Remarks}
\newcommand{\brks}{\begin{rmks} \upshape}
\newcommand{\erks}{\end{rmks}}

\newtheorem{exe}[thm]{Example}
\newcommand{\bexe}{\begin{exe} \upshape}
\newcommand{\eexe}{\end{exe}}

\newtheorem{exes}[thm]{Examples}
\newcommand{\bexes}{\begin{exes} \upshape}
\newcommand{\eexes}{\end{exes}}

\newtheorem*{pre}{Proof}
\newcommand{\bp}{\begin{pre} \upshape}
\newcommand{\ep}{\hfill \qed \end{pre}}
\newcommand{\epp}{\end{pre}}

\newcommand{\beq}{\begin{eqnarray*}}
\newcommand{\eeq}{\end{eqnarray*}}

\newcommand{\beqn}{\begin{equation}}
\newcommand{\eeqn}{\end{equation}}

\newcommand{\ben}{\begin{enumerate}}
\newcommand{\een}{\end{enumerate}}
\newcommand{\bit}{\begin{itemize} \renewcommand{\labelitemi}{$\bullet$} \renewcommand{\labelitemii}{$\star$}}
\newcommand{\eit}{\end{itemize}}

\newcommand{\bfg}{
\begin{figure}[H]
\begin{center}}
\newcommand{\efg}{
\end{center}
\end{figure}
\FloatBarrier}

\newcolumntype{M}[1]{>{\raggedright}m{#1}}

\newcommand{\R}{\mathbb{R}}
\newcommand{\Q}{\mathbb{Q}}
\newcommand{\N}{\mathbb{N}}
\newcommand{\Z}{\mathbb{Z}}

\newcommand{\K}{\mathbb{K}}

\renewcommand{\SS}{\mathbb{S}}

\newcommand{\F}{\mathbb{F}}

\newcommand{\bs}{\symbol{92}}

\newcommand{\ov}{\overline}

\renewcommand{\tilde}{\widetilde}

\renewcommand{\max}{\operatorname{max}}

\newcommand{\rk}{\operatorname{rk}}

\newcommand{\St}{\operatorname{St}}

\newcommand{\eps}{\varepsilon}
\newcommand{\st}{\, | \,}

\newcommand{\ra}{\rightarrow}

\newcommand{\ral}[1]{\underset{#1}{\longrightarrow}}

\newcommand{\liml}{\lim\limits}

\newcommand{\spine}{\operatorname{Spine}}

\newcommand{\f}{\frac}

\renewcommand{\geq}{\geqslant}
\renewcommand{\leq}{\leqslant}

\newcommand{\GL}{\operatorname{GL}}
\newcommand{\SL}{\operatorname{SL}}
\newcommand{\SO}{\operatorname{SO}}

\newcommand{\<}{\langle}
\renewcommand{\>}{\rangle}
\newcommand{\pif}{{+\infty}}

\newcommand{\mk}{\medskip}
\newcommand{\cal}{\mathcal}

\makeatletter
\def\Ddots{\mathinner{\mkern1mu\raise\p@
\vbox{\kern7\p@\hbox{.}}\mkern2mu
\raise4\p@\hbox{.}\mkern2mu\raise7\p@\hbox{.}\mkern1mu}}
\makeatother

\usetikzlibrary{arrows}

\title{A link condition for simplicial complexes, and CUB spaces}
\author{Thomas Haettel}
\date{\today}

\begin{document}

\selectlanguage{english}

\maketitle

\begin{center}
\begin{minipage}{0.8\textwidth}
\textsc{Abstract.}
We advertise the study of metric spaces with a unique convex geodesic bicombing, which we call CUB spaces. These encompass many classical notions of nonpositive curvature, such as CAT(0) spaces and Busemann-convex spaces. Groups having a geometric action on a CUB space enjoy numerous properties.

We want to know when a simplicial complex, endowed with a natural polyhedral metric, is CUB. We establish a link condition, stating essentially that the complex is locally a lattice. This generalizes Gromov's link condition for cube complexes, for the $\ell^\infty$ metric.

The link condition applies to numerous examples, including Euclidean buildings, simplices of groups, Artin complexes of Euclidean Artin groups, (weak) Garside groups, some arcs and curve complexes, and minimal spanning surfaces of knots.
\end{minipage}
\end{center}

\let\thefootnote\relax\footnotetext{Thomas Haettel, thomas.haettel@umontpellier.fr, IMAG, Univ Montpellier, CNRS, France, and IRL 3457, CRM-CNRS, Universit\'{e} de Montr\'{e}al, Canada.\\
{\bf Keywords} : Nonpositive curvature, link condition, simplicial complexes, injective metric spaces, lattices, buildings, braid groups, Artin groups, Garside groups, curve complexes, complexes of groups. {\bf AMS codes}~: 20F65,20F67,20F36,05B35,06A12,20E42,52A35}

\section*{Introduction}

Nonpositive curvature has proven to be a very rich way to study large families of groups. Classical Riemannian sectional nonpositive curvature has been well-studied, and we are interested in developing analogous tools in the setting of cell complexes. More precseily, assume that $X$ is a cell complex, where each cell is identified with a convex polytope in $\R^n$, endowed with some norm. If $X$ is piecewise Euclidean, one may ask whether $X$ is locally CAT(0), which serves as a perfect analogue of nonpositive curvature, and has many examples and applications (\cite{bridson_haefliger}, \cite{bgs}, \cite{ballmann_lectures_NPC}, \cite{caprace_monod_CAT0_structure}, \cite{caprace_monod_cat0_discrete}, \cite{duchesne_CAT0_survey}, \cite{mccammond_CAT0_survey} among many others). However, even if there is a metric link criterion ensuring that a space is locally CAT(0), it is not a combinatorial criterion. The only general situation where such a criterion becomes combinatorial is the setting of CAT(0) cube complexes, i.e. cell complexes made out of unit Euclidean cubes.

\bthn[{Gromov's link condition \cite{gromov_hyperbolic_groups},\cite[Theorem~B.8]{leary_CCC}}]\

A cube complex $X$ is locally CAT(0) if and only if the link of every vertex of $X$ is a flag simplicial complex.
\ethn

The theory of CAT(0) cube complexes, notably due to this extremely simple criterion, knows a huge success, popularized by Agol's and Wise's works (notably) leading to the solution of the virtual Haken conjecture for $3$-manifolds (\cite{haglund_wise}, \cite{haglund_wise_ccc}, \cite{wise_book_hierarchy},\cite{agol}). There are also numerous works on CAT(0) cube complexes (\cite{sageev}, \cite{caprace_sageev}, \cite{chatterji_niblo}, \cite{niblo_roller}, \cite{roller}, \cite{niblo_reeves}, among many others).

\mk

However, outside of the world of cube complexes, it becomes excessively hard to decide whether a given cell complex is CAT(0). For instance, braid groups act properly and cocompactly by isometries on a nice simplicial complex, the dual Garside complex, which is conjectured to be CAT(0). However, this question appears to be quite hard to answer (\cite{brady_mccammond}, \cite{b6}, \cite{jeong}). And one cannot replace this space with a CAT(0) cube complex for braid groups (\cite{haettel_artin_cubic}).

\mk

One problem is that the class of CAT(0) cube complexes, though very rich, remains smaller than the class of spaces we want to call nonpositively curved. For instance, any group with Kazhdan's property (T) may only act with a fixed point on any such CAT(0) cube complex, whereas many of them, such as cocompact lattices in higher rank simple Lie groups, do act properly and cocompactly on a CAT(0) space.

\mk

We are therefore interested in describing criteria ensuring that a simplicial complex, with a norm on each simplex, has nonpositive curvature in some sense. We suggest the following notion of nonpositive curvature: we call a metric space CUB, or Convexly Uniquely Bicombable, if it admits a unique convex geodesic bicombing. It encompasses some other notions of nonpositive curvature, such as CAT(0) and Busemann-convex, and it retains many important properties of classical nonpositive curvature, see Section~\ref{sec:CUB_spaces} for precisions. The study of geodesic bicombings is quite interesting, and has recently developed (see~\cite{wenger},\cite{kasprowski_rueping},\cite{lang}, \cite{descombes_lang_flats}, \cite{descombes_lang_hyperbolicity}, \cite{descombes_asymptotic_rank}, \cite{kleiner_lang_higher_rank_hyperbolicity}, \cite{miesch}, \cite{basso_bicombings} notably).

\mk

Note that one would like to allow such a theory to encompass the case of a simplicial tiling of $\R^n$, for all $n \geq 1$, which is the "zero curvature" case.

One simple combinatorial notion of nonpositive curvature is that of systolic complexes (\cite{januskiewicz_swiatkowski_systolic}, \cite{przytycki_systolic_fixed_point}, \cite{osajda_przytycki_systolic_boundary}); however, one cannot one tile $\R^n$ by a systolic complex, for $n \geq 3$. There is not even a proper action of $\Z^3$ on a systolic complex.

A naive approach would be to consider the Euclidean metric of a standard unit simplex, where each edge has length $1$, and wonder whether the resulting metric is CAT(0). However, it only works well in dimensions $1$ and $2$: indeed, according to~\cite{Kpoczynski_pak_przytycki_acute_triangulations}, any regular CAT(0) simplicial complex is systolic. And if $\R^n$ is tiled by Euclidean simplices with acute angles, then $n \leq 3$.

\mk

These negative results suggest that we need to consider simplicial complexes where simplices have a particular symmetry. We will consider two such symmetries:
\bit
\item The case where simplices have a cyclic order symmetry, which we call type A, which corresponds for instance to the case of the $\tilde{A}_n$ simplex tiling of $\R^n$, see Figure~\ref{fig:A2tilde_pavage}.
\item The case where simplices have a total order symmetry, which we call type C, which corresponds for instance to the case of the $\tilde{C}_n$ simplex tiling of $\R^n$, see Figure~\ref{fig:C2tilde_pavage}.
\eit

In type A, we define a specific polyhedral metric in Section~\ref{sec:norms_on_simplices}, which turns out to behave exceptionally well for applications. In type C, we define the $\ell^\infty$ orthoscheme metric, which is recalled in Section~\ref{sec:norms_on_simplices} and has already been studied in~\cite{brady_mccammond} (for the $\ell^2$ version) and in~\cite{haettel_injective_buildings} and \cite{haettel_helly_kpi1}. This set-up allows us to state our main results in a loose way, just in order to get the general philosophy. We refer the reader to Section~\ref{sec:results_criteria_CUB} for the precise (yet simple) statements of Theorem~\ref{thm:main_criterion_type_A_CUB}, Theorem~\ref{thm:main_criterion_type_C_CUB} and Theorem~\ref{thm:main_criterion_Garside_flag_CUB}.

\bthn[Link condition for simplicial complexes] Let $X$ denote a simplicial complex with totally ordered or cyclically ordered simplices. Then $X$ is locally CUB if and only if $X$ is locally a lattice.
\ethn

Amazingly enough, it turns out that a large number of simplicial complexes appear to satisfy this local lattice property. In particular, this applies to classical Euclidean buildings, to Artin complexes of some Euclidean Artin groups, to some arc complexes of punctured spheres, to Garside and weak Garside groups, to general simplices of groups, to the complex of homologous multicurves on a surface, to the Kakimizu complex of minimal Seifert surfaces of a link. See Section~\ref{sec:applications} for all these applications. Also note that the lattice property is actually simple to verify: one has to check that the poset does not contain any bowtie, see Section~\ref{sec:results_criteria_CUB}.

\mk

Note that this theory applies to the Garside complex of the braid group for either the standard or the dual Garside structure (\cite[Theorem~E]{haettel_helly_kpi1}); whereas with the Euclidean metric, the standard Garside complex is not CAT(0) (\cite[Theorem~3.11]{haettel_helly_kpi1}), and the dual Garside complex is only conjectured to be CAT(0) (\cite{brady_mccammond,b6}).

\mk

The idea of considering spaces which locally are lattices appears in~\cite{brady_mccammond}, \cite{chalopin_chepoi_hirai_osajda}, \cite{b6}, \cite{hirai_CAT0} and \cite{hirai_uniform_modular}. It is also one of the key components of~\cite{haettel_injective_buildings}, \cite{haettel_helly_kpi1} and \cite{haettel_huang_weakly_modular}.

\mk

One key ingredient for the main theorem comes from the theory of injective metric spaces, applied to various spaces constructed from lattices, see~\cite{haettel_injective_buildings} and \cite{haettel_helly_kpi1} which essentially provide the existence of a convex bicombing. The uniqueness part is new, and constitutes the main technical part of this article. It relies notably on work of Descombes and Lang on the combinatorial dimension (\cite{descombes_lang_hyperbolicity}). Interestingly, one of the arguments involves the horofunction boundary of a vector space with a polyhedral norm.

\mk

In addition to the new ideas involved in the proofs of the main theorems, there are several motivations for this article. One is to advertise the study of convex geodesic bicombings on spaces, which looks rich and promising. Another one is to illustrate how results from~\cite{haettel_helly_kpi1}, \cite{descombes_lang_hyperbolicity} and \cite{miesch} can be rendered accessible through simple combinatorial criteria for simplicial complexes. A third motivation is to illustrate the simplicity of these criteria through numerous examples. We believe these combinatorial link conditions will be useful in future works.

\mk

Here is the organization of the article. In Section~\ref{sec:CUB_spaces} and \ref{sec:CUB_groups}, we define CUB spaces and review examples and properties of CUB spaces and groups acting on them. In Section~\ref{sec:norms_on_simplices}, we make precise which polyhedral norms on simplicial complexes we will consider. In Section~\ref{sec:results_criteria_CUB}, we state precisely the various link conditions for locally CUB simplicial complexes. In Section~\ref{sec:lattices_orthoschemes_injective}, we recall definitions of orthoscheme complexes of lattices and their relationship with injective metric spaces. In Section~\ref{sec:bicombings_on_quotient}, we focus on the uniqueness of the convex bicombing in the diagonal quotient of the orthoscheme complex of a graded lattice. In Section~\ref{sec:completion_of_posets}, we show that this result extends to a possibly non-graded lattice. Finally in Section~\ref{sec:proof_local_criteria}, we complete the proofs of the link conditions. In Section~\ref{sec:applications}, we list numerous applications of these link conditions.

\mk

{\bf Acknowledgments:} The author would like to thank Pierre-Emmanuel Caprace, Victor Chepoi, Anthony Genevois, Jingyin Huang, Dan Margalit, Jon McCammond, Piotr Przytycki, Andrew Putman and Constantin Vernicos for very interesting discussions. The author would also like to thank warmly the anonymous referee for his/her very careful reading and his/her numerous insightful comments.

\tableofcontents

\section{Bicombings and CUB spaces} \label{sec:CUB_spaces}

Let $X$ denote a geodesic metric space. Typically $X$ will not be uniquely geodesic, hence the need to select, for each pair of points, a geodesic between them. This is called a geodesic bicombing.

\mk

In this article, a \emph{convex bicombing} on $X$ will be a convex, consistent, reversible geodesic bicombing, i.e. a map $\sigma:X \times X \times [0,1] \ra X$ such that:
\ben
\item (Bicombing) For each $x,y \in X$, the map $t \in [0,1] \mapsto \sigma(x,y,t)$ is a constant speed reparametrized geodesic from $\sigma(x,y,0)=x$ to $\sigma(x,y,1)=y$.
\item (Reversible) For each $x,y \in X$, for each $t \in [0,1]$, we have $\sigma(x,y,t)=\sigma(y,x,1-t)$.
\item (Consistent) For each $x,y \in X$, for each $s,t \in [0,1]$, we have $\sigma(x,y,st)=\sigma(x,\sigma(x,y,t),s)$.
\item (Convex) For each $x,y,x',y' \in X$, the map $t \in [0,1] \mapsto d(\sigma(x,y,t),\sigma(x',y',t)$ is convex (see Figure~\ref{fig:combing}).
\een

\begin{figure}
\begin{center}
\begin{tikzpicture}
\def \p {0.05}
\def \op {1}
\def \gris {black!10}

\draw[fill] (-3,1) circle (\p) node(a) {};
\node (a) at ([yshift=0.5cm]a) {\bfseries $x$};
\draw[fill] (3,1) circle (\p) node(b) {};
\node (b) at ([yshift=0.5cm]b) {\bfseries $y$};
\draw[fill] (-4,-0.5) circle (\p) node(a') {};
\node (a') at ([yshift=-0.5cm]a') {\bfseries $x'$};
\draw[fill] (3.5,-0.5) circle (\p) node(b') {};
\node (b') at ([yshift=-0.5cm]b') {\bfseries $y'$};

\draw [blue, thick] plot [smooth, tension=1] coordinates { (-3,1) (0,0.5) (3,1)};
\draw [blue, thick] plot [smooth, tension=1] coordinates { (-4,-0.5) (0,0) (3.5,-0.5)};

\draw[fill] (0,0.5) circle (\p) node(c) {};
\node (c) at ([yshift=0.5cm]c) {\bfseries $\sigma(x,y,t)$};
\draw[fill] (0,0) circle (\p) node(c') {};
\node (c') at ([yshift=-0.5cm]c') {\bfseries $\sigma(x',y',t)$};

\draw [red, thick] plot [smooth, tension=1] coordinates { (0,0) (0,0.5)};

\end{tikzpicture}
\end{center}
\caption{A convex bicombing}
\label{fig:combing}
\end{figure}

A weaker but useful notion is that of a \emph{conical} bicombing, i.e. satisfying
$$\forall x,x',y,y' \in X, \forall t \in [0,1], d(\sigma(x,y,t),\sigma(x',y',t)) \leq (1-t)d(x,x')+td(y,y').$$
Note that any conical, consistent bicombing is convex.

Moreover, a very important but simple remark is that, if $\sigma$ is a convex (or conical) bicombing on $X$, then any ball in $X$ is $\sigma$-convex.

\mk

If $\sigma$ is a bicombing on $X$, for any $x,y \in X$, we will often denote by $\sigma(x,y)$ the image of the geodesic $\sigma(x,y,\cdot) : [0,1] \ra X$.

\mk

A \emph{local convex bicombing} on $X$ is, for every $x \in X$, the choice of a neighbourhood $U$ of $x$ such that $U$, with the induced metric, has a convex bicombing. Let us recall that, in this article, a neighbourhood of $x$ is any subset of $X$ containing an open subset containing $x$. Since balls are convex with respect to convex bicombings, note that $U$ can be chosen to be a small ball centered at $x$.
\mk

\bdf[CUB]\

We say that a metric space $X$ is \emph{convexly uniquely bicombable} (CUB) if $X$, and each ball of $X$, admit a unique convex bicombing.

We say that a metric space $X$ is \emph{locally convexly uniquely bicombable} (locally CUB) if every point in $X$ admits a neighbourhood which is CUB for the subspace metric.
\edf

Remark that a metric space $X$ is locally CUB if and only if every point in $X$ admits a basis of neighbourhoods consisting of balls which admit a unique convex bicombing.

In the definition of CUB, we need to ask the property for each ball of $X$ in order to ensure that each CUB space is locally CUB. The fact that CUB serves as notion of nonpositive curvature is first justified by the following Cartan-Hadamard type theorem.

\bthm[{\cite{miesch}}] \label{thm:miesch_CUB_cartan_hadamard}
Let $X$ denote a complete, simply-connected, locally CUB metric space. Then $X$ is CUB.
\ethm

Note that Miesch's result is not stated for CUB spaces, but for spaces with a consistent local convex bicombing. The local CUB property is one way to ensure the consistency of a local convex bicombing.

\mk

This notion also encompasses numerous spaces of nonpositive curvature. 

\bexes\
\bit
\item Any CAT(0) space is CUB.
\item A uniquely geodesic metric space is CUB if and only if it is Busemann-convex.
\item A Riemannian manifold is locally CUB if and only if it has nonpositive sectional curvature.
\item A finite-dimensional piecewise Euclidean cell complex, endowed with the length metric, is CUB if and only if it is CAT(0).
\item A finite-dimensional piecewise hyperbolic cell complex, endowed with the length metric, is CUB if and only if it is CAT(-1).
\item Any proper, finite-dimensional injective metric space is CUB (\cite{descombes_lang_hyperbolicity}).
\item Any Banach space is CUB, where the unique convex bicombing is the affine bicombing (\cite[Corollary~1.3]{basso_bicombings}).
\eit
\eexes

Let us comment on two examples. If $X$ is a Riemannian manifold, then it is locally uniquely geodesic. Therefore it is locally CUB if and only if it is Busemann-convex. According to the second variation formula, this is equivalent to having nonpositive sectional curvature (see~\cite{busemann_npc}, and also~\cite{ballmann_lectures_NPC}, \cite{bgs} or \cite{bridson_haefliger} for instance).

If $X$ is a piecewise Euclidean or piecewise hyperbolic cell complex, assume the local CUB property, we want to see that $X$ is CAT(0) or CAT(-1). In particular, $X$ has a continuous geodesic bicombing. Then, by induction on dimension, we can use Bowditch's criterion (\cite{bowditch}, or \cite{bridson_haefliger} in the case $X$ has finitely many shapes) to reduce to proving that the link of every vertex has no locally geodesic loop of length smaller than $2\pi$. The existence of such a loop would give rise to a discontinuity of geodesics, so vertex links are CAT(1), and $X$ is CAT(0) or CAT(-1).

\mk

Let us remind the reader that a geodesic metric space is called \emph{injective} if any family of pairwise intersecting balls has a non-empty global intersection, see for instance~\cite{lang}. It turns out that these metric spaces play a key role underlying the proofs of the main results of this article, relying notably on~\cite{haettel_helly_kpi1}. See Section~\ref{sec:lattices_orthoschemes_injective} for a brief overview.

\mk

Let us state very general properties of CUB metric spaces, which are typical of nonpositive curvature.

\bthm \label{thm:topological_properties_CUB}
Let $X$ denote any CUB metric space. We have the following.
\ben
\item $X$ is contractible.
\item $X$ admits Euclidean isoperimetric inequalities.
\item $X$ admits a Z-boundary and a Tits boundary. Moreover, any action on $X$ by isometries extends continuously to the boundaries.
\een
\ethm

\bp
Let us give a reference for each statement.
\ben
\item Since $X$ admits a convex bicombing, it retracts continuously to any point, and so $X$ is contractible.
\item According to~\cite{wenger}, any metric space with a convex geodesic bicombing admits Euclidean isoperimetric inequalities.
\item According to~\cite[Theorem~1.4]{descombes_lang_hyperbolicity}, any metric space $X$ with a convex geodesic bicombing admits a Z-boundary $\partial X$, defined as equivalence classes of rays with the cone topology. See~\cite{bestvina_zstructure} for the definition of a Z-boundary. If we endow this boundary with the finer Tits topology, we can use results from~\cite{kleiner_lang_higher_rank_hyperbolicity}.
\een
\ep

As we will see in the sequel, CUB spaces behave almost as good as general CAT(0) spaces, notably regarding all arguments using uniqueness and convexity of geodesics. They also happen to be much more frequent than CAT(0) spaces, as the list of examples above suggests.

\section{Groups acting on CUB spaces} \label{sec:CUB_groups}

We now turn to properties of isometric group actions on CUB spaces. First we have a simple fixed point result, due to Basso (\cite{basso_fixed_point}).

\bpro \label{pro:action_CUB_fixed_point}
Let $X$ denote any complete CUB metric space. If $G$ is a group of isometries of $X$ with a relatively compact orbit, then $G$ has a fixed point.

Moreover, if the fixed point set of $G$ is not empty, it is contractible.
\epro

\bp
In case $G$ has a finite orbit, there is a straightforward argument. According to~\cite[Theorem~2.1]{descombes_asymptotic_rank}, there is a well-defined notion of barycenter on $X$, which is equivariant under isometries. In particular, if $G$ is a finite group of isometries, the barycenter of any orbit is a fixed point.

In case $G$ has a relatively compact orbit, the barycenter argument is more subtle, and we need to use~\cite[Theorem~1.1]{basso_fixed_point}: since $X$ is a CUB metric space, the unique convex bicombing $\sigma$ on $X$ satisfies the assumptions that $\sigma$ has the midpoint property and $G$ preserves $\sigma$.

If $G$ acts on $X$ with a fixed point, let $x,y \in X^G$. Then by uniqueness of the convex bicombing $\sigma$ on $X$, we know that $g \cdot \sigma(x,y)=\sigma(x,y)$. In particular $\sigma(x,y)$ is entirely contained in $X^G$, so $\sigma$ restricts to a convex bicombing on $X^G$. In particular, $X^G$ is contractible.
\ep

\brk
Note that the fixed point result does not extend for bounded orbits, like it does for CAT(0) metric spaces. Indeed, according to~\cite[Section~2]{basso_fixed_point}, there exists a bounded, complete, uniquely geodesic, Busemann-convex metric space $X$ (in particular, a CUB space), and an isometry of $X$ without fixed points.
\erk

In particular, we deduce the following.

\bcor
If a group $G$ acts properly by isometries on a complete CUB space $X$, then $X$ is a classifying space for proper actions of $G$.
\ecor

Moreover, in case we have a proper and cocompact action, we can deduce numerous properties for the group.

\bthm
Let $G$ denote a group acting properly and cocompactly by isometries on a CUB space.
\ben
\item $G$ is semihyperbolic in the sense of Alonso-Bridson, and in particular:
\bit
\item Any polycylic subgroup subgroup of $G$ is virtually abelian.
\item Any finitely generated abelian subgroup of $G$ is quasi-isometrically embedded.
\item The word and conjugacy problems are soluble for $G$.
\item The centralizer of any element of $G$ is finitely generated, quasi-isometrically embedded and semihyperbolic.
\eit
\item $G$ has finitely many conjugacy classes of finite subgroups.
\item $G$ satisfies the Farrell-Jones conjecture.
\item $G$ has type $F_\infty$. If $G$ is torsion-free, it has type $F$.
\item $G$ has at most Euclidean Dehn function and higher filling functions (see~\cite{bestvina_eskin_wortman,leuzinger_young}).
\item $G$ satisfies the coarse Baum-Connes conjecture.
\item $G$ has an EZ-boundary.
\item $G$ has contractible asymptotic cones.
\een
\ethm

\bp
Let $X$ denote the CUB space on which $G$ acts properly and cocompactly by isometries. We will give references for all statements.
\ben
\item Since $X$ admits a convex bicombing and $G$ acts properly coboundedly on $X$, we deduce by~\cite[Theorem~III.$\Gamma$.4.7]{bridson_haefliger} that $G$ is semihyperbolic. For the consequences, see~\cite[Theorem~III.$\Gamma$.4.9, Theorem~III.$\Gamma$.4.10, Proposition~III.$\Gamma$.4.15, Proposition~III.$\Gamma$.4.17]{bridson_haefliger}, .
\item According to Proposition~\ref{pro:action_CUB_fixed_point}, any finite subgroup has a fixed point. Since the action of $G$ is proper and cocompact, there exist finitely many conjugacy classes of finite subgroups.
\item According to~\cite{kasprowski_rueping}, any group acting geometrically on a space with a convex bicombing satisfies the Farrell-Jones conjecture.
\item The proofs of~\cite[Proposition~II.5.13, Lemma~I.7A.15]{bridson_haefliger} adapt from CAT(0) spaces to CUB spaces immediately. Hence the quotient $X/G$ has the homotopy type of a finite $CW$-complex.
\item According to Theorem~\ref{thm:topological_properties_CUB}, $X$ has Euclidean isoperimetric inequalities. These translate to the Dehn function and the higher filling functions of $G$ being at most Euclidean.
\item According to~\cite[Theorem~1.3]{fukaya_oguni}, groups acting geometrically on spaces with convex bicombings satisfy the coarse Baum-Connes conjecture.
\item According to Theorem~\ref{thm:topological_properties_CUB}, $X$ admits a Z-boundary $\partial X$, that is $G$-equivariant. This defines an EZ-boundary as in~\cite{farrell_lafont_EZ}.
\item The group $G$ is quasi-isometric to the CUB space $X$. And any asymptotic cone of $X$ admits a bicombing, and is then contractible.
\een
\ep

\bexes\ 
Here is a list of groups acting properly and cocompactly by isometries on a CUB space.
\bit
\item CAT(0) groups, which include fundamental groups of compact nonpositively curved Riemannian manifolds, and uniform lattices in semisimple Lie groups over local fields.
\item Gromov-hyperbolic groups (\cite{lang}), and groups hyperbolic relative to Helly groups (\cite{osajda_valiunas}).
\item Helly groups (\cite{chalopin_chepoi_hirai_osajda}), which include FC type Artin groups and weak Garside groups (\cite{huang_osajda_helly}), and more generally injective groups (\cite{haettel_injective_buildings},\cite{haettel_helly_kpi1}).
\eit
\eexes

Note that mapping class groups of surfaces, and more generally hierarchically hyperbolic groups, satisfy a weaker property: they act properly and coboundedly by isometries on an injective metric space. Such a space admits a (non-unique) conical bicombing, but possibly no convex bicombing, see~\cite{haettel_hoda_petyt}.

\mk

It turns out that many arguments about CAT(0) spaces rely only on the fact that they are uniquely geodesic with associated convex bicombing. These arguments often carry on to the CUB setting, as the following nonpositive criterion for complexes of groups. Let us remind the reader that a scwol is a simple category without loops (see~\cite{bridson_haefliger}). It is quite a general framework to study complexes of groups.

\bthm[{\cite[Theorem~III.C.4.17]{bridson_haefliger}}] \label{thm:complex_groups_CUB_developable}

Let $G({\cal Y})$ denote a complex of groups over a scwol ${\cal Y}$ such that the geometric realization $|{\cal Y}|$ has a metric such that, for each $\sigma \in {\cal Y}$, the local development at $\sigma$ is locally CUB. Then $G({\cal Y})$ is developable, and the simply-connected development of ${\cal Y}$ is CUB.
\ethm

\section{Norms on simplicial complexes} \label{sec:norms_on_simplices}

We are interested in finding nice shapes of simplices that allow to tile $\R^n$, for any $n \geq 1$. We noted in the introduction that the unit Euclidean $n$-simplex does not tile $\R^n$ if $n \geq 3$ (\cite{Kpoczynski_pak_przytycki_acute_triangulations}). This suggests considering simplices with a particular symmetry.

\mk

As a consequence of the classification of Euclidean Coxeter groups (see for instance~\cite[Table~6.1]{davis_coxeter}), there are $4$ infinite families of Euclidean tilings of $\R^n$, for all $n \geq 1$, with fundamental domain a simplex. They correspond to the Coxeter complexes of the Euclidean Coxeter groups of types $\tilde{A}_n$, $\tilde{B}_n$, $\tilde{C}_n$ and $\tilde{D}_n$. Note that both $\tilde{B}_n$ and $\tilde{D}_n$ tilings may be refined into the $\tilde{C_n}$ tiling. See Figures~\ref{fig:A2tilde_pavage} and \ref{fig:C2tilde_pavage} for pictures of the $\tilde{A}_2$ tiling and the $\tilde{C}_2$ tiling.

\begin{figure}[H]
\centering
\begin{tikzpicture}
\def \p {0.05}
\def \op {0.1}
\def \gris {black}
\clip (-3.2,-3.2) rectangle (3.2,3.2);

\draw[fill,color=blue!20!white] (0,0) -- (1,0.58) -- (0,1.16) -- cycle;

\foreach \i in {-10,...,10}
\draw (\i,-10) -- (\i,10);

\foreach \i in {-10,...,10}
\draw (-17.32-\i/2,-10+\i*0.866) -- (17.32-\i/2,10+\i*0.866);

\foreach \i in {-10,...,10}
\draw (-17.32-\i/2,10-\i*0.866) -- (17.32-\i/2,-10-\i*0.866);

\end{tikzpicture}
\caption{The $\tilde{A}_2$ tiling of the plane, with the $\tilde{A}_2$ simplex in blue.}
\label{fig:A2tilde_pavage}
\end{figure}

\begin{figure}[H]
\centering
\begin{tikzpicture}
\def \p {0.05}
\def \op {0.1}
\def \gris {black}
\clip (-3.2,-3.2) rectangle (3.2,3.2);

\draw[fill,color=blue!20!white] (0,0) -- (1,1) -- (1,0) -- cycle;

\foreach \i in {-10,...,10}
\draw (\i,-10) -- (\i,10);

\foreach \i in {-10,...,10}
\draw (-10,\i) -- (10,\i);

\foreach \i in {-10,...,10}
\draw (-10,-10+2*\i) -- (10,10+2*\i);

\foreach \i in {-10,...,10}
\draw (10,-10+2*\i) -- (-10,10+2*\i);

\end{tikzpicture}
\caption{The $\tilde{C}_2$ tiling of the plane, with the $\tilde{C}_2$ simplex in blue.}
\label{fig:C2tilde_pavage}
\end{figure}

\mk

We will therefore focus on one hand on the simplices coming from the $\tilde{A}_n$ tiling, whose vertices have a total cyclic order, which we will call type A. On the other hand, we will focus on the simplices coming from the $\tilde{C}_n$ tiling, whose vertices have a total order, which we will call type C. We give details below, starting with the type C case which is a bit simpler.

\mk

We will see in the applications (see Section~\ref{sec:applications}) that the assumption that simplices have an order on their vertices is often quite natural and geometric. Oftentimes, these orders will come from a naturally defined rank or type of vertices.

\subsection{Type C}

\begin{figure}[H]
\centering
\begin{tikzpicture}
\def \p {0.05}
\def \op {0.3}
\def \gris {blue}
\draw[fill] (0,0) circle (\p) node(0) {};
\draw[fill] (3,0) circle (\p) node(1) {};
\draw[fill] (3,3) circle (\p) node(2) {};
\draw[fill] (3,3) + (20:3) circle (\p) node(3) {};

\draw[black,fill opacity=\op,fill=\gris] (0.center) -- (1.center) -- (2.center) -- (0.center);
\draw[black,fill opacity=\op,fill=\gris] (1.center) -- (2.center) -- (3.center) -- (1.center);
\draw[dashed] (0.center) -- (3.center);
\draw[-{Stealth[scale=2]}] (0.center) -- (1.center);
\draw[-{Stealth[scale=2]}] (1.center) -- (2.center);
\draw[-{Stealth[scale=2]}] (2.center) -- (3.center);

\node at ([yshift=-0.5cm]0) {\bfseries $v_0=(0,0,0)$};
\node at ([yshift=-0.5cm]1) {\bfseries $v_1=(1,0,0)$};
\node at ([yshift=0.7cm]2) {\bfseries $v_2=(1,1,0)$};
\node at ([yshift=0.5cm]3) {\bfseries $v_3=(1,1,1)$};

\end{tikzpicture}
\caption{The standard $3$-orthosimplex, with the total order on vertices.}
\label{fig:3_orthoscheme}
\end{figure}
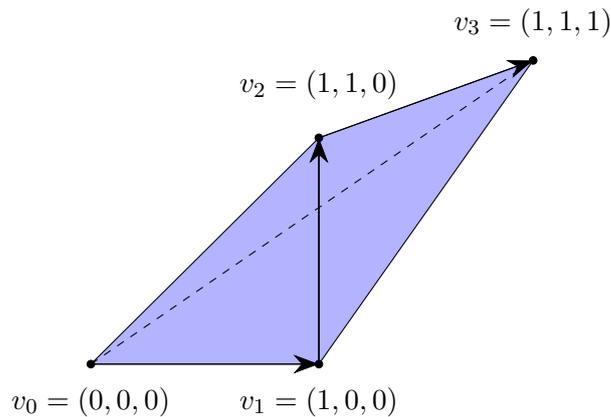

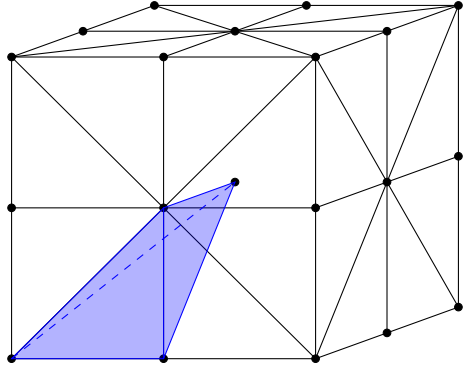
\begin{figure}[H]
\centering
\begin{tikzpicture}[scale = 1]
\def \p {0.05}
\def \r {2}
\def \m {1}
\def \a {20}
\def \op {0.5}
\def \gris {black!10}
\draw[fill] (0,0) circle (\p) node(000) {};
\draw[fill] (4,0) circle (\p) node(100) {};
\draw[fill] (0,4) circle (\p) node(010) {};
\draw[fill] (4,4) circle (\p) node(110) {};
\draw[fill] (2,0) circle (\p) node(200) {};
\draw[fill] (0,2) circle (\p) node(020) {};
\draw[fill] (2,4) circle (\p) node(210) {};
\draw[fill] (4,2) circle (\p) node(120) {};
\draw[fill] (2,2) circle (\p) node(220) {};
\draw[fill] (2,2) + (\a:\m) circle (\p) node(222) {};
\draw[fill] (4,0) + (\a:\r) circle (\p) node(101) {};
\draw[fill] (0,4) + (\a:\r) circle (\p) node(011) {};
\draw[fill] (4,4) + (\a:\r) circle (\p) node(111) {};
\draw[fill] (4,0) + (\a:\m) circle (\p) node(102) {};
\draw[fill] (4,2) + (\a:\m) circle (\p) node(122) {};
\draw[fill] (4,2) + (\a:\r) circle (\p) node(121) {};
\draw[fill] (4,4) + (\a:\m) circle (\p) node(112) {};
\draw[fill] (0,4) + (\a:\m) circle (\p) node(012) {};
\draw[fill] (2,4) + (\a:\r) circle (\p) node(211) {};
\draw[fill] (2,4) + (\a:\m) circle (\p) node(212) {};

\draw (000.center) -- (100.center) -- (110.center) -- (010.center) -- (000.center);
\draw (000.center) -- (110.center);
\draw (100.center) -- (010.center);
\draw (200.center) -- (210.center);
\draw (020.center) -- (120.center);

\draw (100.center) -- (101.center) -- (111.center) -- (110.center);
\draw (100.center) -- (111.center) -- (010.center);
\draw (101.center) -- (110.center) -- (011.center);
\draw (010.center) -- (011.center) -- (111.center);
\draw (012.center) -- (112.center) -- (102.center);
\draw (210.center) -- (211.center);
\draw (120.center) -- (121.center);

\draw[blue, dashed] (000.center) -- (222.center);
\draw[blue,fill opacity=0.3,fill=blue] (220.center) -- (000.center) -- (200.center) -- (220.center) -- (222.center) -- (200.center);


\end{tikzpicture}
\caption{The partition of a cube in $\R^3$ into standard orthosimplices.}
\label{fig:partition}
\end{figure}

Consider an $n$-dimensional simplex $\sigma$, with a total order on its vertices, which we could then label $v_0<v_1< \dots <v_n$. Then one can naturally identify $\sigma$ with the standard $n$-simplex of type $\tilde{C}_n$, which is also called the standard orthosimplex. It may be defined as the convex hull in $\R^n$ of the set of points $v_0=(0,0,\dots,0),v_1=(1,0,\dots,0),\dots,v_n=(1,1,\dots,1)$, see Figure~\ref{fig:3_orthoscheme}. It also coincides with a simplex of the barycentric subdivision of the $n$-cube $[0,2]^n$, see Figure~\ref{fig:partition}.

\mk

As a subset of $\R^n$, the standard orthosimplex may naturally be endowed with one of the following two norms:
\bit
\item The standard $\ell^2$ Euclidean norm.
\item The standard $\ell^\infty$ norm, given by $\|x\|_\infty = \max(|x_1|,|x_2|,\dots,|x_n|)$.
\eit

While many works have focused on the $\ell^2$ Euclidean norm on ortoscheme complexes (see~\cite{brady_mccammond}, \cite{b6}, \cite{jeong}, \cite{dougherty_mccammond_witzel}), in this article we will only consider the standard $\ell^\infty$ norm. Let us be precise about the metric simplicial complexes we will be considering.

\bdf
A simplicial complex $X$ is said to have \emph{ordered simplices} if each simplex of $X$ has a total order on its vertex set, which is consistent with respect to inclusions of simplices.
\edf

For such complexes, we are able to define the standard $\ell^\infty$ metric.

\bdf
Let $X$ denote a simplicial complex with ordered simplices, with finite simplices. Endow each $d$-simplex of $X$ with the standard $\ell^\infty$ norm of the standard orthosimplex in $\R^d$, and endow $X$ with the associated length metric: this is called the \emph{standard $\ell^\infty$ metric}. 
\edf

A very important remark is that this metric is well-defined: indeed if $\tau$ is a face of a simplex $\sigma$ with totally ordered vertices, then the standard $\ell^\infty$ metric on $\tau$ (with the induced order) is an isometric subspace of the standard $\ell^\infty$ metric on $\sigma$. Note that this would not be the case with the Euclidean metric: if $\sigma$ is a triangle with vertices $v_0<v_1<v_2$, then the edge between $v_0$ and $v_2$ has length $\sqrt{2}$ for the standard Euclidean metric of $\sigma$. For the $\ell^\infty$ metric, all edges have length $1$.

\mk

Also note that if $X$ is a simplicial complex with ordered simplices, and if $G$ is a group of simplicial automorphisms which either preserve or reverse the orders on vertices of simplices, then $G$ acts by isometries on $X$ with the standard $\ell^\infty$ metric.

\subsection{Type A}

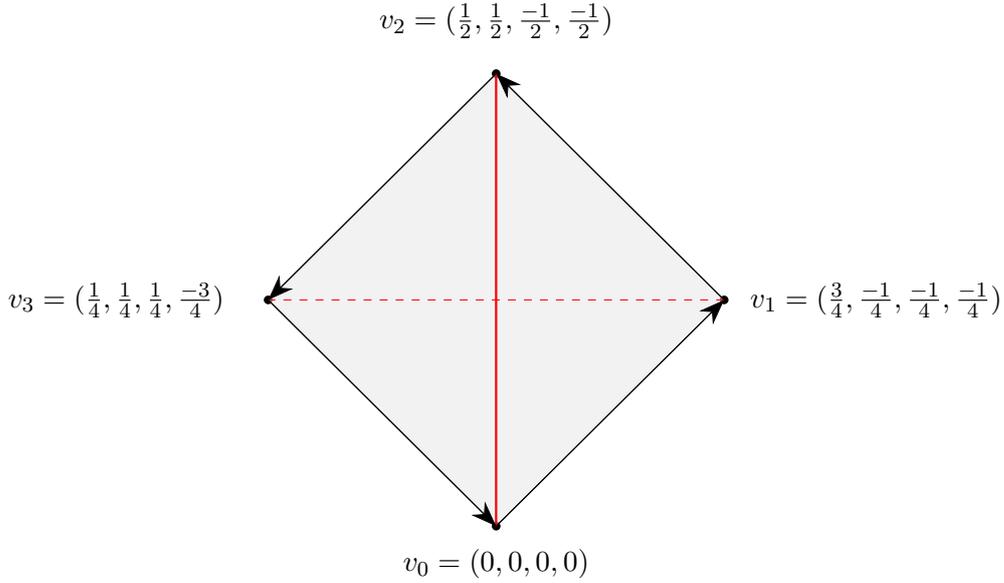
\begin{figure}[H]
\centering
\begin{tikzpicture}
\def \p {0.05}
\def \op {0.5}
\def \gris {black!10}
\draw[fill] (0,0) circle (\p) node(0) {};
\draw[fill] (3,3) circle (\p) node(1) {};
\draw[fill] (0,6) circle (\p) node(2) {};
\draw[fill] (-3,3) circle (\p) node(3) {};

\draw[black,fill opacity=\op,fill=\gris] (0.center) -- (1.center) -- (2.center) -- (0.center);
\draw[black,fill opacity=\op,fill=\gris] (2.center) -- (3.center) -- (0.center) -- (2.center);
\draw[red,thick] (0.center) -- (2.center);
\draw[red,dashed] (1.center) -- (3.center);
\draw[-{Stealth[scale=2]}] (0.center) -- (1.center);
\draw[-{Stealth[scale=2]}] (1.center) -- (2.center);
\draw[-{Stealth[scale=2]}] (2.center) -- (3.center);
\draw[-{Stealth[scale=2]}] (3.center) -- (0.center);

\node at ([yshift=-0.5cm]0) {\bfseries $v_0=(0,0,0,0)$};
\node at ([xshift=2cm]1) {\bfseries $v_1=(\f{3}{4},\f{-1}{4},\f{-1}{4},\f{-1}{4})$};
\node at ([yshift=0.7cm]2) {\bfseries $v_2=(\f{1}{2},\f{1}{2},\f{-1}{2},\f{-1}{2})$};
\node at ([xshift=-2cm]3) {\bfseries $v_3=(\f{1}{4},\f{1}{4},\f{1}{4},\f{-3}{4})$};

\end{tikzpicture}
\caption{The standard $3$-simplex of type $\tilde{A}_3$, with the cyclic order on vertices. The two red edges have dihedral angle $\f{\pi}{2}$, the other four edges have dihedral angle $\f{\pi}{3}$.}
\label{fig:3_orthosimplex_A3}
\end{figure}

Consider an $n$-dimensional simplex $\sigma$, with a total cyclic order on its vertices, which we could then label $(v_i)_{i \in \Z/n\Z}$. Then one can naturally identify $\sigma$ with the standard $n$-simplex of type $\tilde{A}_n$. It may be defined as the convex hull in $\R^n=\{x \in \R^{n+1} \st x_1+x_2+\dots+x_{n+1}=0\}$ of the set of points
$$v_0=(0,0,\dots,0),v_1=(\f{n}{n+1},\f{-1}{n+1},\dots,\f{-1}{n+1}),\dots,v_n=(\f{1}{n+1},\dots,\f{1}{n+1},\f{-n}{n+1}).$$
More precisely, for $0 \leq i \leq n$, let $j=n+1-i$, then the first $i$ coordinates of $v_i$ are equal to $\frac{j}{n+1}$ and the last $j$ coordinates of $v_i$ are equal to $\frac{-i}{n+1}$, see Figure~\ref{fig:3_orthosimplex_A3}.

\mk

In this description, the cyclic symmetry is not transparent: the standard simplex of type $\tilde{A}_n$ may equivalently be defined as
$$\sigma=\{x \in \R^{n+1} \st x_1+x_2+\dots+x_{n+1}=0, x_1 \geq x_2 \geq \dots \geq x_n \geq x_{n+1} \geq x_1-1\}.$$
Therefore, we can also describe the standard simplex of type $\tilde{A}_n$ as the image, in the quotient under the diagonal translation action of $\R$, of the standard column $C$ in $\R^{n-1}$:
$$C=\{x \in \R^{n+1} \st x_1 \geq x_2 \geq \dots \geq x_n \geq x_{n+1} \geq x_1-1\}.$$
Note that the column $C$ has a natural description as a simplicial complex, with linearly ordered $(n+1)$-dimensional simplices with totally ordered vertices. Roughly speaking, these simplices appear when looking at chambers obtained by cutting $C$ using the hyperplanes $\{x_i=k\}_{1 \leq i \leq n+1,k \in \Z}$ in $\R^{n+1}$.

\begin{figure}[H]
\centering
\begin{tikzpicture}
\def \p {0.05}
\def \op {0.5}
\def \gris {black!10}
\def \r {2}
\def \o {20}

\foreach \i in {0,...,3}
\draw[fill] (\r*\i,\r*\i) + (20:\r*\i) circle (\p) node(A\i) {};

\foreach \i in {0,...,3}
\draw[fill] (\r*\i+\r,\r*\i) + (20:\r*\i) circle (\p) node(B\i) {};

\foreach \i in {0,...,3}
\draw[fill] (\r*\i+\r,\r*\i+\r) + (20:\r*\i) circle (\p) node(C\i) {};

\draw[black,fill opacity=\op,fill=\gris] (A0.center) -- (B0.center) -- (C0.center) -- (A0.center);

\draw[black,fill opacity=\op,fill=\gris] (B0.center) -- (C0.center) -- (B1.center) -- (B0.center);
\draw[black,fill opacity=\op,fill=\gris] (B1.center) -- (C1.center) -- (B2.center) -- (B1.center);
\draw[black,fill opacity=\op,fill=\gris] (B2.center) -- (C2.center) -- (B3.center) -- (B2.center);

\draw[black,fill opacity=\op,fill=\gris] (C0.center) -- (B1.center) -- (C1.center) -- (C0.center);
\draw[black,fill opacity=\op,fill=\gris] (C1.center) -- (B2.center) -- (C2.center) -- (C1.center);
\draw[black,fill opacity=\op,fill=\gris] (C2.center) -- (B3.center) -- (C3.center) -- (C2.center);

\foreach \i in {0,...,3}{
\draw[-{Stealth[scale=2]},dashed] (A\i.center) -- (B\i.center);
\draw[-{Stealth[scale=2]}] (B\i.center) -- (C\i.center);}

\foreach \i in {1,...,3}
\draw[dashed] (A\i.center) -- (C\i.center);

\draw[dashed] (A0.center) -- (A1.center) -- (A2.center) -- (A3.center);
\draw[dashed] (B0.center) -- (A1.center);
\draw[dashed] (B1.center) -- (A2.center);
\draw[dashed] (B2.center) -- (A3.center);

\draw[-{Stealth[scale=2]},dashed] (C0.center) -- (A1.center);
\draw[-{Stealth[scale=2]},dashed] (C1.center) -- (A2.center);
\draw[-{Stealth[scale=2]},dashed] (C2.center) -- (A3.center);

\end{tikzpicture}
\caption{The standard column $C$ in $\R^3$, tiled by standard $3$-orthosimplices.}
\label{fig:column_C_R3}
\end{figure}
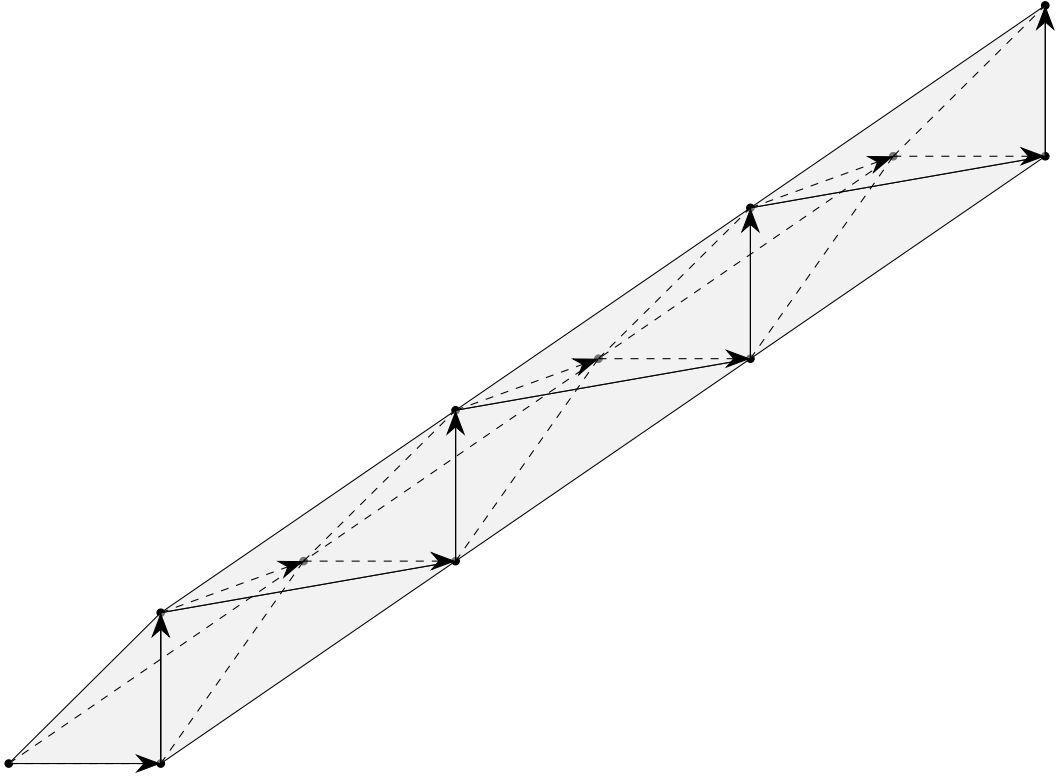

More precisely, we have $C=\bigcup_{q \in \Z, 0 \leq r \leq n} \tau_{q,r}$, where
$$\tau_{q,r} =\{x \in \R^{n+1} \st q \geq x_r \geq x_{r+1} \geq x_{r+2} \geq \dots \geq x_{n+1} \geq x_1-1 \geq x_2-1 \geq \dots \geq x_r-1 \geq q-1\},$$
see Figure~\ref{fig:column_C_R3}. The study and use of columns is fundamental in~\cite{brady_mccammond}, \cite{dougherty_mccammond_witzel} and \cite{haettel_helly_kpi1}.

\mk

It turns out that there are two natural norms on the standard $\tilde{A}_n$ simplex:
\bit
\item The standard $\ell^2$ Euclidean norm.
\item The norm $\|x\| = \sup_{1 \leq i \neq j \leq n+1} |x_i-x_j|$, which we will call \emph{standard polyhedral norm}.
\eit

Note that the standard polyhedral norm on $\R^n$ can also be described more naturally as the quotient of $\R^{n+1}$, endowed with the standard $\ell^\infty$ metric, under the diagonal translation action by $\R$. The unit ball of the standard polyhedral norm is also the quotient of the standard column $C \subset \R^{n+1}$, under the diagonal translation action by $\R$. Note that the standard Euclidean metric on $\R^n$ is also the quotient of the standard Euclidean metric on $\R^{n+1}$.

\mk

Also note that both metrics are invariant under a (possibly order-reversing) cyclic permutation of the vertices of the standard $\tilde{A}_n$ simplex. This is easier to see in the column picture, since for instance $(x_1,x_2,\dots,x_{n+1}) \in C \mapsto (x_2,x_3,\dots,x_{n+1},x_1-1) \in C$ is an isometry for both the standard Euclidean and the standard polyhedral norm. However, an order-reversing permutation of the vertices of the standard $\tilde{A}_n$ simplex does not lift to an isometry of the standard column $C$, but it still induces an isometry of the standard $\tilde{A}_n$ simplex for either norm.

\mk

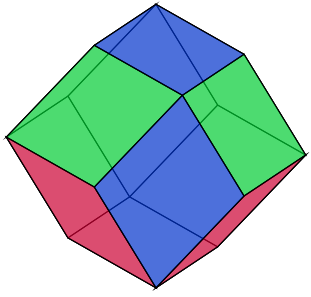
\begin{figure}[H]
\centering
\tdplotsetmaincoords{70}{10}
\begin{tikzpicture}[tdplot_main_coords]

\def \r {2}

\foreach \x in {-1,1}{
\foreach \y in {-1,1}{
\foreach \z in {-1,1}
\coordinate (\x\y\z) at (\x,\y,\z);}}

\foreach \x in {-1,1}{
\coordinate (X\x) at (\x * \r,0,0);
\coordinate (Y\x) at (0,\x * \r,0);
\coordinate (Z\x) at (0,0,\x * \r);}

\foreach \x in {-1,1}{
\foreach \y in {-1,1}{
\draw (X\x) -- (\x\y1) -- (Y\y) -- (\x\y-1) -- cycle;
\draw (X\x) -- (\x1\y) -- (Z\y) -- (\x-1\y) -- cycle;
\draw (Y\x) -- (1\x\y) -- (Z\y) -- (-1\x\y) -- cycle;
}}

\draw [fill opacity=0.7,fill=green!80!blue] (X1) -- (111) -- (Y1) -- (11-1) -- cycle;
\draw [fill opacity=0.7,fill=green!80!blue] (X-1) -- (-111) -- (Y1) -- (-11-1) -- cycle;
\draw [fill opacity=0.7,fill=red!80!blue] (X1) -- (11-1) -- (Z-1) -- (1-1-1) -- cycle;
\draw [fill opacity=0.7,fill=red!80!blue] (X-1) -- (-11-1) -- (Z-1) -- (-1-1-1) -- cycle;
\draw [fill opacity=0.7,fill=blue!80!green] (Y1) -- (111) -- (Z1) -- (-111) -- cycle;
\draw [fill opacity=0.7,fill=blue!80!green] (Y1) -- (11-1) -- (Z-1) -- (-11-1) -- cycle;

\end{tikzpicture}
\caption{The unit ball of the standard polyhedral norm in $\R^3$: the rhombic dodecahedron.}
\label{fig:rhombic_dodecahedron}
\end{figure}

The unit ball of the standard polyhedral norm in $\R^n$ is the projection of the standard $(n+1)$-cube (with edge lengths $2$) in $\R^{n+1}$ by the diagonal translation action of $\R$. In dimension $n=2$, it is a regular hexagon, and in dimension $n=3$ it is called the rhombic dodecahedron, see Figure~\ref{fig:rhombic_dodecahedron}.

\mk

One should also note that the $n$-simplex, endowed with the Hilbert metric (\cite{nussbaum_hilbert_metric}, \cite{lemmens_walsh_isometries_hilbert}), is isometric to the normed space
$$ \left(\R^n=\{x \in \R^{n+1} \st x_1+x_2+\dots+x_{n+1}=0\}, \|x\| = \sup_{1 \leq i \neq j \leq n+1} |x_i-x_j|\right).$$
So the norm is the same as the standard polyhedral norm we defined. However, the embedding of an $n$-simplex in this model is not immediate: in particular, it depends on the particularly chosen cyclic order on the vertices on the $n$-simplex.

\mk

We now make precise how to define a metric if there is a cyclic order on the vertices of each simplex of a simplicial complex.

\bdf
A simplicial complex $X$ is said to have \emph{cyclically ordered simplices} if each simplex of $X$ has a total cyclic order on its vertex set, which is consistent with respect to inclusions of simplices.
\edf

For such complexes, we are able to define the standard polyhedral metric.

\bdf
Let $X$ denote a simplicial complex with cyclically ordered simplices, with finite simplices. Endow each $d$-simplex of $X$ with the standard polyhedral norm of the standard $d$-simplex, and endow $X$ with the associated length metric: this is called the \emph{standard polyhedral metric}. 
\edf

Note that, as in the type C case, this metric is well-defined: if $\tau$ is a face of a simplex $\sigma$ with cyclically ordered vertices, then the standard polyhedral metric on $\tau$ (with the induced cyclic order) is an isometric subspace of the standard polyhedral metric on $\sigma$.

\mk

Also note that if $X$ is a simplicial complex with cyclically ordered simplices, and if $G$ is a group of simplicial automorphisms which either preserve or reverse the cyclic orders on vertices of simplices, then $G$ acts by isometries on $X$ with the standard polyhedral metric.

\section{Statements of the link conditions} \label{sec:results_criteria_CUB}

We will present the link conditions for nonpositive curvature. These are local criteria ensuring the CUB property in the case of simplicial complexes with ordered or cyclically ordered simplices. In order to state the criteria precisely, we will recall basic definitions on lattices and bowties.

\subsection{Lattices and bowties}

We start by recalling necessary definitions of posets, lattices and bowties.

\bdf\

A \emph{lattice} is a poset $L$ such that any two $x,y \in L$ have a minimal upper bound denoted $x \vee y$, called the \emph{join} of $x$ and $y$, and a maximal lower bound denoted $x \wedge y$, called the \emph{meet} of $x$ and $y$.

A meet-\emph{semilattice} is a poset where one only requires the existence of meets.

A join-\emph{semilattice} is a poset where one only requires the existence of joins.

A poset $L$ is \emph{bounded} if it has a minimum (usually denoted $0$), and a maximum (usually denoted $1$).

In a poset $L$, if $x,y \in L$, we say that $y$ \emph{covers} $x$ if $x<y$, and there exists no $z \in L$ with $x<z<y$.

A poset $L$ is \emph{graded} if there exists a \emph{rank} function $\rk : L \ra \Z$ such that, for any $x,y \in L$ such that $y$ covers $x$, we have $\rk(y)=\rk(x)+1$.

A poset $L$ is \emph{homogeneous} if, given any $x<y$ in $L$, there is a bound on the length of chains from $x$ to $y$ (depending on $x,y$).

A graded lattice $L$ has rank $n$ if every maximal chain in $L$ has $n+1$ elements.

A \emph{bowtie} in a poset $L$ consists of $4$ elements $a,b,c,d \in L$ such that $a,b < c,d$, and such that there exists no $x \in L$ such that $a,b \leq x \leq c,d$. If $L$ is graded, we say that a bowtie $a,b < c,d$ is \emph{balanced} if $\rk(a)=\rk(b)$ and $\rk(c)=\rk(d)$.
\edf

\bexe\
Here are simple examples of such posets.
\bit
\item The boolean poset ${\cal P}(\{1,\dots,n\})$, ordered by inclusion, is a bounded graded lattice of rank $n$.
\item The poset of vector subspaces of an $n$-dimensional vector space, ordered by inclusion, is a bounded graded lattice of rank $n$.
\item The poset of finite-dimensional vector subspaces of an arbitrary vector space, ordered by inclusion, is a graded meet-semilattice.
\item The poset of partitions of $\{1,\dots,n\}$, ordered by refinement, is a bounded graded lattice of rank $n-1$.
\item Fix $n \geq 2$, and consider the vertex set $U_n \subset \R^2$ of a regular $n$-gon in the plane. Say that a partition $P$ of $U_n$ is noncrossing if, for any distinct $A,B \in P$, the convex hulls of $A$ and $B$ do not intersect. Then the poset of noncrossing partitions of $U_n$, ordered by refinement, is a bounded graded lattice of rank $n-1$. Note that his example generalizes to any finite Coxeter group (\cite{brady_watt}, \cite{bessis}).
\eit
\eexe

It turns out that bowties are quite efficient to decide whether a poset is a lattice. If $L$ is any poset, let us denote by $L \cup \{0,1\}$ the bounded poset obtained by adding two new elements: a minimum $0$, and a maximum $1$.

\bpro \label{pro:lattice_balanced_bowtie}
Let $L$ denote a graded poset. Then $L \cup \{0,1\}$ is a lattice if and only if $L$ has no balanced bowtie.
\epro

\bp
Assume that $L \cup \{0,1\}$ is a lattice, and consider $a,b<c,d$ in $L$. Then the meet $x$ of $c,d$ is such that $a,b \leq x \leq c,d$. So $L$ has no bowties.

\mk

Conversely, assume that $L$ has no bowtie $a,b < c,d$ with $\rk(a)=\rk(b)$ and $\rk(c)=\rk(d)$.

Fix $c,d \in L$ arbitrary with $\rk(c)=\rk(d)$, we will prove that $c,d$ have a meet in $L \cup \{0\}$. It is sufficient to prove that no bowtie in $L$ contains $c,d$ as upper elements. For a contradiction, assume that there exists a bowtie $a,b < c,d$. By assumption, we know that $\rk(a) \neq \rk(b)$, say $\rk(a) < \rk(b)$. Pick $b' \in L$ such that $b' < b$ and $\rk(b')=\rk(a)$. By assumption, there exists $x \in L$ such that $a,b' \leq x \leq c,d$. Since $a,b<c,d$ is a bowtie, we deduce that $a=x$. So $a \leq b$, which contradicts that $a,b<c,d$ is a bowtie.

\mk

By symmetry, we also know that any $a,b \in L$ with $\rk(a)=\rk(b)$, have a join in $L \cup \{1\}$.

\mk

Now fix any $c,d \in L$, we will prove that $c,d$ have a meet in $L \cup \{0\}$. As before, it is sufficient to prove that no bowtie in $L$ contain $c,d$ as upper elements. For a contradiction, assume that there exists a bowtie $a,b < c,d$.

If $\rk(a)=\rk(b)$, then since $a,b$ have a join this contradicts that $a,b<c,d$ is a bowtie.

So we know that $\rk(a) \neq \rk(b)$, say $\rk(a) < \rk(b)$. Pick $b' \in L$ such that $b' < b$ and $\rk(b')=\rk(a)$. Now let $x$ denote the join of $a,b'$: it is such that $a,b' \leq x \leq c,d$. Since $a,b<c,d$ is a bowtie, we deduce that $a=x$. So $a \leq b$, which contradicts that $a,b<c,d$ is a bowtie.

Hence any $c,d \in L$ have a meet in $L \cup \{0\}$.

\mk

By symmetry, we also know that any $a,b \in L$ have a join in $L \cup \{1\}$. Hence $L \cup \{0,1\}$ is a lattice.
\ep

\subsection{Type $A$}

Given a flag simplicial complex $X$ with cyclically ordered simplices, for each vertex $x \in X$, there is an antisymmetric relation $\leq_x$ on $\St(x)$, with minimal element $x$, obtained by declaring that $y \leq_x z$ if and only if $x,y,z$ form a triangle, and $x,y,z$ are cyclically ordered with respect to the original cyclic ordering on $X$.

\bdf
We say that a simplicial complex with cyclically ordered simplices is a \emph{local poset} if, for any vertex $x \in X$, the pair $(\St(x),\leq_x)$ is a poset. In other words:
$$\forall x,y,z,t \in X, (x \leq y \leq z \leq x) \wedge (x \leq z \leq t \leq x) \Rightarrow (x \leq y \leq t \leq x).$$ 
\edf

We are now able to state precisely the link condition for simplicial complexes with cyclically ordered simplices. Let us remind that a cell complex $X$ is called \emph{locally finite-dimensional} if, for any cell $C$ of $X$, there is a bound on the dimension of cells containing $C$.
 
\bmthm \label{thm:main_criterion_type_A_CUB}
Let $X$ denote a locally finite-dimensional flag simplicial complex with cyclically ordered simplices, which is a local poset. Then $X$, endowed with the standard polyhedral norm, is locally CUB if and only if for every vertex $x \in X$, the poset $(\St(x),\leq_x)$ is a meet-semilattice.
\emthm

Note that, in case $(\St(x),\leq_x)$ is not a meet-semilattice, we prove that $X$ does not even have one local conical bicombing.

\mk

Let us remark that the semilattice condition in Theorem~\ref{thm:main_criterion_type_A_CUB} has several equivalent formulations.

\blem \label{lem:characterizations_lattice_type_A}
Let $X$ denote a locally finite-dimensional flag simplicial complex with cyclically ordered simplices, which is a local poset. For each vertex $x \in X$, the following are equivalent.
\ben
\item $(\St(x),\leq_x)$ is a meet-semilattice.
\item $(\St(x) \cup \{1\},\leq_x)$ is a lattice.
\item $(\St(x), \leq_x)$ has no bowtie.
\een
In the case where all maximal simplices containing $x$ have the same dimension, the poset $(\St(x), \leq_x)$ is graded, and these conditions are also equivalent to:
\ben
\setcounter{enumi}{3}
\item $(\St(x), \leq_x)$ has no balanced bowtie.
\een
\elem

\bp\
Let us denote by $L$ the poset $(\St(x),\leq_x)$. We will first prove the equivalence between properties (1),(2) and (3).

\bit
\item [$(1) \Rightarrow (2)$] Let us assume that $L$ is a meet-semilattice, and consider $y,z \in L$. Since there is a bound on length of chains in $L$, let us assume that $a,b \in L \cup \{1\}$ are minimal among common upper bounds of $y,z$ in $L \cup \{1\}$. If $a \neq b$, then $c=a \wedge b$ is an upper bound for $y,z$ such that $c < a,b$, which is a contradiction. Hence $a=b$ is the join of $y$ and $z$. So $L \cup \{1\}$ is a lattice.
\item [$(2) \Rightarrow (3)$] Let us assume that $L \cup \{1\}$ is a lattice. If $a,b<c,d$ is a bowtie in $L$, then $y=c \wedge d$ is such that $a,b \leq y \leq c,d$, which contradicts the bowtie property. So $L$ has no bowtie.
\item [$(3) \Rightarrow (1)$] Let us assume that $L$ has no bowtie. As in the proof of $(1) \Rightarrow (2)$, let $y,z \in L$, and consider $a,b \in L$ maximal among common lower bounds of $y,z$. Since $a,b \leq y,z$ and $L$ has no bowtie, we deduce that there exists $c \in L$ such that $a,b \leq x \leq y,z$. By maximality of $a,b$, we deduce that $a=b=c$, which is the meet of $y,z$. Hence $L$ is a meet-semilattice.
\eit

Assume now that all maximal simplices containing $x$ have the same dimension $n$. Then the poset $L$ is graded and has rank $n$. According to Proposition~\ref{pro:lattice_balanced_bowtie}, properties (2) and (4) are equivalent.
\ep

\subsection{Type $C$}

Given a flag simplicial complex $X$ with ordered simplices, for each vertex $x \in X$, there is a natural induced antisymmetric relation $\leq_x$ on $\St(x)$.

\bdf
We say that a simplicial complex with ordered simplices is a \emph{local poset} if, for any vertex $x \in X$, the set $(\St(x),\leq_x)$ is a poset. In other words:
$$\forall x \in X, \forall y,z,t \in \St(x), (y \leq z) \wedge (z \leq t) \Rightarrow (y \leq t).$$ 
\edf

Given a simplicial complex $X$ with ordered simplices, for each vertex $x \in X$, there is an antisymmetric relation $\leq_x$ on $\St(x)$, obtained by declaring that $y \leq_x z$ iff $x,y,z$ form a triangle, and the edge between $y$ and $z$ is oriented from $y$ to $z$ with respect to the original ordering on $X$.

We will also consider the ascending star $\St^+(x)=\{y \in \St(x) \st y \geq x\}$, with the induced order $\leq_x$, whose minimum is $x$. And the descending star $\St^-(x)=\{y \in \St(x) \st y \leq x\}$, with the induced order $\leq_x$, whose maximum is $x$.

\mk

We are now able to state precisely the link condition for simplicial complexes with ordered simplices.

\bmthm \label{thm:main_criterion_type_C_CUB}
Let $X$ denote a locally finite-dimensional flag simplicial complex with ordered simplices, which is a local poset. Then $X$, endowed with the standard $\ell^\infty$ norm, is locally CUB and is locally injective if and only if for any vertex $x \in X$, we have:
\bit
\item (Lattice condition) The poset $(\St(x),\leq_x)$ has no bowtie.
\item (Flag condition) Any $a,b,c \in \St(x)$ which are pairwise upperly bounded (resp. lowerly bounded) have a common upper bound (resp. lower bound).
\eit
\emthm

Remark that, in the case $(\St(x),\leq_x)$ does not satisfy any of the two conditions, we prove that $X$ is not locally injective and does not even have one local convex bicombing. Note that this flag condition is called "lattice-theoretic flag condition" in~\cite{hirai_CAT0}.

\mk

As in the cyclically ordered case, let us remark that the lattice condition has several equivalent formulations. The proof is similar to Lemma~\ref{lem:characterizations_lattice_type_A}.

\blem \label{lem:characterizations_lattice_type_C}
Let $X$ denote a locally finite-dimensional flag simplicial complex with ordered simplices, which is a local poset. For each vertex $x \in X$, the following are equivalent.
\ben
\item $(\St(x), \leq_x)$ has no bowtie.
\item $(\St^+(x) ,\leq_x)$ is a meet-semilattice and $(\St^-(x) ,\leq_x)$ is a join-semilattice.
\item $(\St(x) \cup \{0,1\},\leq_x)$ is a lattice.
\een

In the case where all maximal simplices containing $x$ have the same dimension, the poset $(\St(x), \leq_x)$ is graded, and these conditions are also equivalent to:
\ben
\setcounter{enumi}{3}
\item $(\St(x), \leq_x)$ has no balanced bowtie.
\item $(\St^+(x) ,\leq_x)$ and $(\St^-(x) ,\leq_x)$ have no balanced bowtie.
\een
Similarly, the flag condition is equivalent to asking that the following two conditions hold:
\bit
\item Any $a,b,c \in \St^+(x)$ which are pairwise upperly bounded have a common upper bound.
\item Any $a,b,c \in \St^-(x)$ which are pairwise lowerly bounded have a common lower bound.
\eit
\elem

Also note that Theorem~\ref{thm:main_criterion_type_C_CUB} applies to all cube complexes. Indeed if $X$ is a locally finite-dimensional cube complex, consider the barycentric subdivision $X'$ of $X$. It is a simplicial complex with ordered simplices, and the standard $\ell^p$ metric on orthosimplices coincides with the standard $\ell^p$ metric on cubes (up to a factor $2$). If $x$ is a vertex of $X$, the lattice condition for $\St_{X'}(x)$ is equivalent to requiring that the link of $x$ in $X$ is simplicial, and the flag condition for $\St_{X'}(x)$ is equivalent to requiring that the link of $x$ in $X$ is a flag simplicial complex. Hence we recover Gromov's link condition for the $\ell^\infty$ norm. Note that, for a finite-dimensional cube complex, requiring that the piecewise Euclidean metric is CAT(0) is equivalent to requiring that the piecewise $\ell^\infty$ metric is CUB (see~\cite{miesch_CCC}).

\mk

The analogous statement for the piecewise Euclidean metric is false, see~\cite[Theorem~3.10]{haettel_helly_kpi1}. However, according to~\cite{hirai_CAT0} (see also~\cite{b6}), it is true under the extra (very restrictive) assumption that $(\St(e),\leq_e)$ is a semimodular lattice.

\subsection{Garside flag complexes}

A particular class of examples encompassing both types of shapes of simplices are Garside flag complexes (see~\cite{haettel_huang_weakly_modular}), which we now define.

\bdf \label{def:garside_flag_complex}
A \emph{Garside flag complex} is a pair $(X,\varphi)$, where $X$ is a simply-connected flag simplicial complex with finite simplices, with ordered simplices, and $\varphi$ is an order-preserving automorphism of $X$, such that the following holds:
\bit
\item For any simplex $\sigma$ of $X$, we have that $\sigma \cup \varphi(\min \sigma)$ is a simplex of $X$.
\item For any vertex $x \in X$, we have $\varphi(x)>x$, and the interval $[x,\varphi(x)]$ is a homogeneous lattice.
\eit
\edf

If $(X,\varphi)$ is a Garside flag complex, the quotient $X/\varphi$ is defined as the flag simplicial complex with vertex set $X^{(0)} / \<\varphi\>$, whose $k$-simplices correspond to images of chains $x_0 < x_1 < \dots < x_k < \varphi(x_0)$. Note that the quotient $X/\varphi$ has cyclically ordered simplices.

\mk

For this class of simplicial complexes, we are now able to state precisely the link condition. 

\bmthm \label{thm:main_criterion_Garside_flag_CUB}
Let $(X,\varphi)$ denote a Garside flag complex. Then $X$, endowed with the standard $\ell^\infty$ metric, is CUB and injective. Moreover the quotient $X/\varphi$, endowed with the standard polyhedral metric, is CUB.
\emthm

Note that the lattice property is also necessary, as in the proof of Theorem~\ref{thm:main_criterion_type_C_CUB}.

\section{Lattices, orthoscheme complexes and injective metric spaces} \label{sec:lattices_orthoschemes_injective}

We will review here the relationship between lattices, orthoscheme complexes and injective metric spaces developed in \cite{haettel_injective_buildings} and \cite{haettel_helly_kpi1}.

\subsection{Injective metric spaces and combinatorial dimension}

A geodesic metric space is called \emph{injective} if any family of pairwise intersecting balls has a non-empty global intersection. See for example~\cite{lang} for an introduction to injective metric spaces. It turns out that injective metric spaces are ubiquitous:

\bthm[\cite{isbell}]
Any metric space $X$ embeds isometrically in a unique minimal injective metric space $EX$, its injective hull.
\ethm

Injective metric spaces are relevant for CUB spaces because Lang proved that any injective metric space admits a canonical conical geodesic bicombing. Under properness assumption, Descombes and Lang improved the result to an actual convex geodesic biombing.

\bthm[{\cite[Theorem~1.1]{descombes_lang_hyperbolicity}}]
Let $X$ denote a proper injective metric space. Then $X$ admits a convex geodesic bicombing.
\ethm

Concerning uniqueness, Descombes and Lang provided a criterion relying on the notion of combinatorial dimension.

\bdf
The \emph{combinatorial dimension} of a metric space $X$ is the topological dimension of its injective hull $EX$.
\edf

\bexe
For instance, if $X$ is a CAT(0) cube complex with the piecewise $\ell^\infty$ metric, the combinatorial dimension of $X$ coincides with its dimension as a cube complex.
\eexe

Note that, if $X$ is an isometric subspace of $Y$, the combinatorial dimension of $X$ is bounded above by the combinatorial dimension of $Y$.

\mk

The combinatorial dimension of a metric space is usually hard to compute. However, we have the following criterion due to Dress.

\bthm[{\cite[Theorem~4.1]{descombes_lang_hyperbolicity}}]
Let $X$ denote a metric space, and let $n \geq 1$ be an integer. The space $X$ has combinatorial dimension at most $n$ if and only if for every finite subset $Z \subset X$ with $2n+2$ elements and for every fixed-point-free involution $i:Z \ra Z$, there exists a fixed-point-free-bijection $j:Z \ra Z$ distinct from $i$ such that
$$\sum_{z \in Z} d(z,i(z)) \leq \sum_{z \in Z} d(z,j(z)).$$
\ethm

Descombes and Lang proved that, for metric spaces with finite combinatorial dimension, there could be at most one convex geodesic bicombing.

\bthm \cite[Theorem~1.2]{descombes_lang_hyperbolicity} \label{thm:descombes_lang_finite_dimension_unique_bicombing}
Let $X$ denote a metric space with finite combinatorial dimension. Then $X$ admits at most one convex geodesic bicombing.
\ethm

Combining both results, we immediately get the following.

\bcor[Descombes-Lang]
Let $X$ denote a proper, finite-dimensional, injective metric space. Then $X$ is CUB.
\ecor

We also record, for later use, the following elementary decomposition result.

\bpro \label{pro:orthoscheme_complex_local_product}
Let $L$ denote a poset, such that $x$ is comparable to every element of $L$. Let $L^+=\{y \in L \st y \geq x\}$ and $L^-=\{y \in L \st y \leq x\}$. Then the geometric realization $|L|$ of $L$, with the standard $\ell^\infty$ metric, is locally isometric at $x$ to the $\ell^\infty$ product $|L^+| \times |L^-|$ of the geometric realizations of $L^+$ and $L^-$ with the standard $\ell^\infty$ metrics.
\epro

\bp
It suffices to notice that this statement is true for orthosimplices, with the $\ell^\infty$ metric. Consider the standard $n$-orthosimplex $\sigma \subset \R^n$, with vertices $v_0=(0,0,\dots,0)<v_1=(1,0,\dots,0)< \dots <v_n=(1,\dots,1)$. More precisely, $\sigma$ is defined by the following inequalities:
$$\sigma=\{x \in \R^n \st 1 \geq x_1 \geq x_2 \geq \dots \geq x_n \geq 0\}.$$ 
Given any $0 < k < n$, note that $v_k$ does not lie on the hyperplane $\{x_k=x_{k+1}\}$ supporting a face of $\sigma$. So a neighbourhood of $v_k$ in $\sigma$ is isometric to a neighourhood of $v_k$ in
$$E=\{x \in \R^n \st 1 \geq x_1 \geq x_2 \geq \dots \geq x_k \geq 0, 1 \geq x_{k+1} \geq x_{k+2} \geq \dots \geq x_n \geq 0\}.$$ 
The space $E$, with the $\ell^\infty$ norm, is isometric to the $\ell^\infty$ product of two standard orthosimplices of dimensions $k$ and $n-k$.
\ep

\subsection{Orthoscheme complexes of lattices}

We now explain that many examples of injective metric spaces come from geometric realizations of posets.

\mk

Let $L$ denote a poset. Then the geometric realization $X$ of $L$ is the simplicial complex whose simplices are chains in $L$. Note that each simplex of $X$ has thus an induced total order on its set of vertices, so it can be endowed with the standard $\ell^\infty$ norm. We will then endow $X$ with the induced length metric.

\mk

One of the main interests in lattices lies in the following result. 

\bthm[{\cite[Theorem~3.9]{haettel_helly_kpi1} \label{thm:orthoscheme_lattice_injective}}]
Let $L$ denote a bounded graded lattice. Then the geometric realization of $L$, with the standard $\ell^\infty$ norm, is injective and CUB.
\ethm

In the sequel, we will need more details about the result in Theorem~\ref{thm:orthoscheme_lattice_injective}. If $L$ is a bounded, graded lattice of rank $n$, we start by describing what we call the affine version of $L$.

\mk

We will define a new poset $M$, which will be called the \emph{affine version} of $L$. Let $C(L)$ denote the set of maximal chains $c_{0,1}=0 <_L c_{1,2} <_L \dots c_{n-1,n} <_L c_{n,n+1}=1$ in $L$.  With these notations, the element denoted $c_{i,i+1}$ has rank $i$.

Let us consider the following subspace of $\R^n$:
$$\Sigma = \{u \in \R^n \st u_1 \leq u_2 \leq \dots \leq u_n\}.$$

\mk

For each maximal chain $c \in C(L)$, let $\Sigma_c$ denote a copy of $\Sigma$.

Let us consider the space
$$M = \bigcup_{c \in C(L)} \Sigma_c / \sim,$$
where for each $c,c' \in C(L)$, if we denote $I=\{1 \leq i \leq n-1 \st c_{i,i+1} \neq c'_{i,i+1}\}$, we identify $\Sigma_c$ and $\Sigma_{c'}$ along the subspaces
$$\{u \in \Sigma_c \st \forall i \in I, u_i=u_{i+1}\} \simeq \{u \in \Sigma'_c \st \forall i \in I, u_i = u_{i+1}\}.$$
We can describe $M$ as a quotient of the space $M_0 = C(L) \times \Sigma$. If $c \in C(L)$ and $u \in \Sigma$, let us denote by $[c,u]$ the equivalence class of $(c,u) \in M_0$ in $M$

\bexe
One illustrating example is the following: consider the boolean lattice $L$ of rank $n$, i.e. the lattice of subsets of the finite set $\{1,\dots,n\}$, with the inclusion order. 
Maximal chains in $L$ correspond to permutations of $\{1,\dots,n\}$. The space $M$ may be identified with $\R^n$, where for each permutation $w$ of $\{1,\dots,n\}$, the subspace $\Sigma_w$ is
$$\Sigma_w=\{x \in \R^n \st x_{w(1)} \leq x_{w(2)} \leq \dots \leq x_{w(n)}\}.$$
\eexe

We will endow $M$ with the length metric induced by the standard $\ell^\infty$ metric on each subspace $\Sigma_c \subset \R^n$, for $c \in C(L)$. Note that the geometric realization $X$ of $L$ is naturally a subspace of $M$:
$$X = \{[c,u] \in M \st c \in C(L), u \in \Sigma, 0 \leq u_1 \leq u_2 \leq \dots u_n \leq 1\}.$$

\mk

This metric space $M$ satisfies the following. Recall that a subset $X$ of a metric space $M$ with a bicombing $\sigma$ is called $\sigma$-convex if, given $x,y \in X$, we have $\sigma(x,y) \subset X$.

\bthm[{\cite[Theorem~3.8]{haettel_helly_kpi1} \label{thm:affine_version_injective}}]
The metric space $M$ is injective, and has a unique convex bicombing $\sigma$. Moreover, the subspace $X \subset M$ is isometric and $\sigma$-convex. 
\ethm

\mk

There is also a diagonal isometric action of $\R$ on $M$, given by $s \cdot [c,u]=[c,u_1+s,u_2+s,\dots,u_n+s]$. We also know that the bicombing on $M$ behaves nicely with respect to this action, as stated below.

\bpro \label{pro:bicombing_translation_invariant}
The convex bicombing $\sigma$ on $M$ satisfies the following:
$$\forall x,x' \in M, \forall t \in [0,1], \forall s,s' \in \R, \sigma(s \cdot x,s' \cdot x',t) = ((1-t)s+ts') \cdot \sigma(x,x',t).$$
\epro

\bp
This is a consequence of the proof of~\cite[Theorem~3.8]{haettel_helly_kpi1}, since the midpoint-nonincreasing conical geodesic bicombing defined there satisfies this equality.
\ep

We also show that stars of simplices are convex in $X$.

\blem \label{lem:stars_convex}
Let $L$ denote a bounded graded lattice, and let $X$ denote the geometric realization of $L$, with the standard $\ell^\infty$ norm. Let $c$ denote a chain in $L$ containing $0$ and $1$, and let $A \subset X$ denote the corresponding simplex of $X$. Then the star of $A$ in $X$ is $\sigma$-convex.
\elem

\bp
For each $a \in c \bs \{0,1\}$, let $A_a \subset A \subset X$ denote the corresponding $2$-simplex of $X$. Note that the star of $A$ in $X$ equals the intersection of the stars of $A_a$, for all $a \in c \bs \{0,1\}$. Therefore it is sufficient to consider the case $c=(0<a<1)$.

\mk

Fix $x,x' \in \St(A)$. By considering $\St(A) \subset X \subset M$, we may find $s \geq 0$ such that $s \cdot x \in \St(A)$ and $s \cdot x \geq a$. Similarly, we may find $s' \geq 0$ such that $s' \cdot x' \in \St(A)$ and $s' \cdot x' \geq a$. Note that the interval $\{x \in X \st x \geq a\}$ is $\sigma$-convex in $X$. So we deduce that, for any $t \in [0,1]$, we have $\sigma(s \cdot x,s' \cdot x',t) \geq a$, so in particular $\sigma(s \cdot x,s' \cdot x',t) \in \St(A)$. According to Proposition~\ref{pro:bicombing_translation_invariant}, we deduce that $\sigma(x,x',t) \in \R \cdot \St (A) \cap X = \St(A)$. Hence $\St(A)$ is $\sigma$-convex.
\ep

\subsection{Unique bicombings in orthoscheme complexes of semilattices}

We now turn to the case of orthoscheme complexes of semilattices. Here, in order to obtain the injectivity and the CUB property, we need to ask for the flag condition.

\bthm \label{thm:orthoscheme_semilattice_injective}
Let $L$ denote a graded meet-semilattice with minimum $0$. Assume that any $a,b,c \in L$ pairwise upperly bounded have a common upper bound. Then the orthoscheme complex of $L$, with the standard $\ell^\infty$ metric, is injective and CUB.
\ethm

\bp
According to~\cite[Theorem~6.1]{haettel_helly_kpi1}, we know that the orthoscheme complex $|L|$ of $L$, with the standard $\ell^\infty$ metric, is injective. Since it is finite-dimensional, according to Theorem~\ref{thm:descombes_lang_finite_dimension_unique_bicombing} we know that it admits at most one convex bicombing. It remains to prove that it admits a convex bicombing.

\mk

Let $L'=L \cup \{1\}$ denote the bounded poset obtained by adding a maximal element $1$. Let $\pi:|L'| \ra |L|$ denote the simplicial map obtained by sending each $x \in L$ to itself, and sending $1$ to $0$. It is a $1$-Lipschitz retraction.

Note that $L'$ is a bounded graded lattice, so according to Theorem~\ref{thm:orthoscheme_lattice_injective} we know that its geometric realization $|L'|$, endowed with the standard $\ell^\infty$ metric, is injective and admits a unique convex bicombing $\sigma'$.

Let us define $\sigma : L \times L \times [0,1] \ra L$ by $\sigma(x,y,t) = \pi(\sigma'(x,y,t))$: since $\pi$ is a $1$-Lipschitz retraction, we deduce that $\sigma$ is geodesic, and that $\sigma$ is conical. Furthermore, according to Proposition~\ref{pro:bicombing_translation_invariant}, we know that $\sigma$ is consistent, hence $\sigma$ is convex. In conclusion, $\sigma$ is a convex bicombing on $|L|$. So $|L|$ is CUB.
\ep

\section{Bicombings on quotients of orthoscheme complexes} \label{sec:bicombings_on_quotient}

We will now explain how unique convex bicombings on orthoscheme complexes give rise to unique convex bicombings on their diagonal quotient. The key result we will use is Descombes and Lang's Theorem~\ref{thm:descombes_lang_finite_dimension_unique_bicombing}. However, bounding the combinatorial dimension is quite subtle.

\subsection{Horofunctions and combinatorial dimension}\ \label{subsec:horofunctions}

\mk

We start by developing an interesting relationship between the horoboundary of a metric space and its combinatorial dimension, which could possibly be used in different contexts. Let us first remind the definition of horoboundary. The idea of representing points of a metric space as distance functions is due to Gromov, and has then been quite studied from various perspectives (see~\cite{borel_ji}, \cite{guivarch}, \cite{brill_horofunction_polyhedral}, \cite{ji_schilling_horofunction_polyhedral}, \cite{haettel_schilling_wienhard}, \cite{walsh_horofunction_normed_space}, \cite{walsh_horofunction_hilbert}, \cite{ciobotaru_kramer_schwer_polyehdral_I}).

Recall that if $X$ is a proper metric space, its \emph{horoboundary} $\partial X$ is defined as the boundary of the image of the embedding
\beq X & \mapsto & \R^X/\R\\
x & \mapsto & d(x,\cdot)+\R,\eeq
where $\R$ acts on $\R^X$ by postcomposition by the standard addition. Note that, if $\xi \in \partial X$, then for any $x,y \in X$, the quantity $\xi(x)-\xi(y)$ is well-defined.

\mk

For each geodesic ray $c:[0,\infty) \ra X$, there is an associated horofunction $\xi_c \in \partial X$ such that
$$\forall x,y \in X, \liml_{t \ra \infty} d(x,c(t))-d(y,c(t)) = \xi_c(x)-\xi_c(y).$$
Note that these horofunctions are particular cases of Busemann points (\cite{walsh_horofunction_normed_space}), where we allow for almost geodesic rays.

\mk

When we restrict to a finite-dimensional real vector space with a polyhedral norm, we have the following.

\bthm[{\cite{brill_horofunction_polyhedral,walsh_horofunction_normed_space,ji_schilling_horofunction_polyhedral}}] \label{thm:horofunction_compactification_vector_space}
Let $X$ denote a finite-dimensional real vector space, with a polyhedral norm. Then its horofunction compactification $X \cup \partial X$ is naturally homeomorphic to the unit ball in the dual space $X^*$.
\ethm

We immediately deduce the following.

\bpro \label{pro:vector_space_finite_Busemann_points}
Let $X$ denote a finite-dimensional real vector space, with a polyhedral norm, and fix $x_0 \in X$. The affine half-lines issued from $x_0 \in X$ correspond to finitely many horofunctions in $\partial X$.
\epro

\bp
According to Theorem~\ref{thm:horofunction_compactification_vector_space}, there are only finitely many horofunctions in $\partial X$, up to translation. And each of these translation classes has a unique representative given by an affine half-line starting from a fixed $x_0 \in X$.
\ep

We now show that a particular space with a finiteness assumption on the horoboundary has finite combinatorial dimension. Let us recall that a bicombing $\sigma$ on $X$ is called geodesically complete if, for any distinct $x,y \in X$, there exists a bi-infinite $\sigma$-line through $x$ and $y$.

\bthm \label{thm:conical_bicombing_combinatorial_dimension_Busemann}
Assume that $X$ is a proper metric space with a geodesically complete convex bicombing $\sigma$. Assume that there exists $x_0 \in X$ such that there are at most $N \in \N$ horofunctions in $\partial X$ corresponding to $\sigma$-rays from $x_0$. Then $X$ has combinatorial dimension at most $N$.
\ethm

We will in fact need a slightly more precise version of Theorem~\ref{thm:conical_bicombing_combinatorial_dimension_Busemann}, which we state below.

\blem \label{lem:conical_bicombing_combinatorial_dimension_Busemann}
Assume that $X$ is a metric space with a convex bicombing $\sigma$. Assume that there exist $x_0 \in X$ and subsets $Y,Z \subset X$, with $x_0 \in Y \subset Z \subset X$, such that the following hold:
\ben
\item $Z$ is proper.
\item For any $x,y \in Y$, there exists an infinite $\sigma$-ray from $x$ through $y$ eventually contained in $Z$.
\item For any $x \in Z$, $\sigma(x_0,x) \subset Z$.
\item There are at most $N \in \N$ horofunctions in $\partial X$ corresponding to $\sigma$-rays from $x_0$ contained in $Z$.
\een
Then $Y$, with the subspace metric, has combinatorial dimension at most $N$.
\elem

\bp
Let $B$ denote the finite set of horofunctions corresponding to $\sigma$-rays from $x_0$. Fix $x,y \in Y$ distinct, we will prove that there exists $\beta \in B$ such that $d(x,y) = D = |\beta(x)-\beta(y)|$. We will then explain why this bounds the combinatorial dimension of $Y$ by $|B|$.

\mk

Consider an infinite $\sigma$-geodesic ray $c:[0,+\infty) \ra X$ starting from $c(0)=x$, with $c(D)=y$, eventually contained in $Z$. Since $Z$ is proper, there exists a sequence $(t_n)_{n \in \N}$ going to $+\infty$ such that $\sigma(x_0,c(t_n))$ converges to an infinite $\sigma$-geodesic ray $c_0:[0,+\infty) \ra Z$ starting from $c_0(0)=x_0$ to some $\beta \in B$.

We know that the two rays $c,c_0$ are at Hausdorff distance at most $C \geq 0$. For each $n \in \N$ large enough, let $z_n \in c_0$ such that $d(x,z_n)=t_n$: we have $d(z_n,x)-d(z_n,y) \ral{n \ra \pif} \beta(x)-\beta(y)$. For each $n \in \N$, let $y_n=\sigma(x,z_n,D)$. Since $d(c(t_n),z_n) \leq 2C$, we deduce that $d(y,y_n) \ral{n \ra \pif} 0$ by conicality. As a consequence, we know that
$$\liml_{n \ra \pif} d(z_n,x)-d(z_n,y) = \liml_{n \ra \pif} d(z_n,x)-d(z_n,y_n) = D = \beta(x)-\beta(y).$$

\mk

We now explain how this helps to bound the combinatorial dimension of $Y$. Fix representatives $B_0$ of $B$. Consider the map
\beq \varphi : Y &\ra& \R^{B_0}\\
x & \mapsto & (\beta \in B_0 \mapsto \beta(x)).\eeq
If we endow $\R^{B_0}$ with the standard $\ell^\infty$ distance, we just proved that the map $\varphi$ is an isometric embedding. Indeed every horofunction is $1$-Lipschitz. Since $\R^{B_0}$ is an injective metric space with dimension $|B|$, we deduce that $Y$ has combinatorial dimension at most $|B|$.
\ep

\subsection{Cone complexes and geodesically complete bicombings}\ \label{subsec:cone_complexes}

\mk

In order to use Theorem~\ref{thm:conical_bicombing_combinatorial_dimension_Busemann}, we need to focus on geodesically complete bicombings. Since we want a local criterion, we will use a framework which allows to replace the neighbourhood of a point in a simplicial complex by the neighbourhood of a point in a cone complex, for which the geodesic completeness is easier to prove. Let us now define more precisely cone polyhedral complexes, as we did in~\cite{haettel_hoda_petyt_Lp_metrics}.

\mk

A \emph{cone polyhedron} $P$ is a finite intersection
$\bigcap_{i=1}^k H_i$ of linear half-spaces of some
finite-dimensional normed space $\bigl(\R^n, |\cdot|\bigr)$ such that
$0 \in \partial H_i$ for all $i$.  A \emph{face} or \emph{cell} of $P$
is a proper subset $P \cap \bigcap_{i \in I} \partial H_i$ for some
$I \subseteq \{1,2, \ldots, k\}$.  The \emph{spine} of $P$ is the face
$P \cap \bigcap_{i=1}^k \partial H_i$, it is also the largest affine
subspace of $P$ that contains $0$.  The spine of any face of $P$
coincides with the spine of $P$.

A \emph{cone polyhedral complex}, or just \emph{cone complex}, with
spine a normed vector space $\spine(X)$, is a metric space obtained, 
from a collection of cone polyhedra whose spines contain an
isometric copy of $\spine(X)$, by gluing along isometries of
their faces that identify $\spine(X)$. The metric on a cone complex is
the length metric induced by its polyhedra.

Fix $x_0 \in X$, let $C$ be the minimal closed cell containing $x_0$
and let $X_C$ be the closed star of $C$.  Given $x \in X_C$ there is a
unique affine map $\phi_x$ from $[0,1]$ to the minimal closed cell of
$X$ that contains both $x_0$ and $x$ such that $\phi_x(0) = x_0$ and
$\phi_x(1) = x$. Let $rx$ denote the image $\phi_x(r)$.  For
$n \geq 1$, let $nX_C$ be $X_C$ with the metric scaled up by $n$.
Then for $n < m$, there is an isometric embedding $nX_C \to mX_C$
given by $x \mapsto \frac{n}{m}x$.  Viewing this isometric embedding
as an inclusion, we define the \emph{tangent cone} $T_{x_0} X$ of $X$
at $x_0$ as the ascending union $\bigcup_{n \geq 1} nX_C$.

We say that a cell complex or a cone complex $X$ is \emph{locally conical} if, for every point $x \in X$, if we denote by $C_x$ the minimal closed cell containing $x$, the open star of $C_x$ is a neighbourhood of $x$ in $X$. For instance, if $X$ has finitely many shapes, then $X$ is locally conical (recall that a cell complex is said to have finitely many shapes if it has finitely many isometry types of cells). 

An immediate consequence of the definition of a locally conical cone complex is the following.

\blem \label{lem:cone_complex_tangent_cone}
Let $X$ denote a locally conical cell complex. Given $x \in X$, let $C_x$ be the minimal closed cell containing $x$. The metric tangent cone $T_x X$ of $X$ at $x$ is a locally conical cone complex such that any small enough ball centered at $x$ in $X$ is isometric to a ball centered at $x$ in $T_x X$. Moreover, the spine of $T_x X$ has the same dimension as $C_x$.
\elem

We can use this definition to study geodesically complete bicombings.

\blem \label{lem:bicombing_finite_number_chambers}
Let $X$ denote a finite-dimensional normed locally conical cone complex, and let $\sigma$ denote a convex bicombing on $X$. Then, for any $x,y \in X$, $\sigma(x,y)$ is contained in a finite union of cells of $X$.
\elem

\bp
Fix $x,y \in X$. For each $z \in \sigma(x,y)$, by local conicality, there exists $\eps>0$ such that the ball $B(z,\eps)$ is a simplicial cone. Therefore we know that $\sigma(x,y) \cap B(z,\eps) = \sigma(x,z) \cup \sigma(z,y)$ is contained in a union of two cells of $X$. The statement follows from compactness of $\sigma(x,y)$.
\ep

Before stating the next result, let us remind that a (constant speed) geodesic $\gamma:[0,1] \ra X$ in a metric space is called \emph{straight} if, for any $z \in X$, the function $t \in [0,1] \mapsto d(\gamma(t),z)$ is convex. For instance, if $\sigma$ is a convex bicombing on $X$, any $\sigma$-line is straight. This has been studied by Descombes and Lang, who proved the following result, which is the main argument in the proof of Theorem~\ref{thm:descombes_lang_finite_dimension_unique_bicombing}, see~\cite[Proposition~4.3]{descombes_lang_hyperbolicity}.

\bpro \label{pro:finite_combinatorial_at_most_one_straight}
Let $X$ denote a metric space of finite combinatorial dimension. Between any two points of $X$, there exists at most one straight geodesic.
\epro

The following consequence is immediate.

\blem \label{lem:finite_dimension_implies_unique_straight_geodesics}
Let $X$ denote a finite-dimensional normed locally conical cone complex. Assume that $X$ admits a convex bicombing $\sigma$. Assume that any finite subcomplex of $X$, with the subspace metric, has finite combinatorial dimension. Then $X$ is CUB and $\sigma$ is the unique convex bicombing on $X$.
\elem

\bp
Assume for a contradiction that $\sigma'$ is another convex bicombing on $X$. Let $x,y \in X$ such that $\sigma(x,y) \neq \sigma'(x,y)$. According to Lemma~\ref{lem:bicombing_finite_number_chambers}, $\sigma(x,y)$ and $\sigma'(x,y)$ are contained in finitely many cells: let $Y$ denote a finite subcomplex of $X$ containing both $\sigma(x,y)$ and $\sigma'(x,y)$. By assumption, the subspace $Y$ has finite combinatorial dimension. According to Proposition~\ref{pro:finite_combinatorial_at_most_one_straight}, there exists at most one straight geodesic in $Y$ between $x$ and $y$. This contradicts that $\sigma(x,y)$ and $\sigma'(x,y)$ are distinct straight geodesics in $Y$ between $x$ and $y$. Hence $\sigma$ is the unique convex bicombing on $X$, and $X$ is CUB.
\ep

We now explain how to bound the combinatorial dimension of a finite subcomplex.

\bthm \label{thm:geod_complete_implies_CUB}
Let $X$ denote a finite-dimensional polyhedrally normed locally conical cone complex. Assume that $X$ admits a convex bicombing $\sigma$. Assume that, for any cells $C,C'$ of $X$, there exists a finite subcomplex $Z$ of $X$ such that, for any $x \in C$ and $y \in C'$, there exists an infinite $\sigma$-ray from $x$ through $y$ eventually contained in $Z$. Then any finite subcomplex of $X$, with the subspace metric, has finite combinatorial dimension. Moreover, $X$ is CUB and $\sigma$ is the unique convex bicombing on $X$.
\ethm

\bp
Let $Y$ denote a finite subcomplex of $X$. By assumption, there exists a finite subcomplex $Z$ of $X$ such that, for any $x,y \in Y$, there exists an infinite $\sigma$-ray from $x$ through $y$ eventually contained in $Z$. Let us endow $Y$ and $Z$ with the subspace metric of $X$.

\mk

Now fix $x_0$ in the spine of $X$. According to Proposition~\ref{pro:vector_space_finite_Busemann_points}, for each cell $C$ of $X$, there are only finitely many horofunctions in $\partial X$ which correspond to $\sigma$-rays from $x_0$ which are contained in $C$. As a consequence, since $Z$ is finite, there are only finitely many horofunctions in $\partial X$ which correspond to $\sigma$-rays from $x_0$ which are contained in $Z$.
\mk

According to Lemma~\ref{lem:conical_bicombing_combinatorial_dimension_Busemann}, we deduce that $Y$ has finite combinatorial dimension. According to Lemma~\ref{lem:finite_dimension_implies_unique_straight_geodesics}, we deduce that $X$ is CUB and $\sigma$ is the unique convex bicombing on $X$.
\ep

We extract from the previous results simple criteria ensuring that a cell complex has finite combinatorial dimension or is CUB. Let us recall that a face in a finite-dimensional complex is \emph{free} if it is properly contained in a unique maximal face.

\bcor
Let $X$ denote a finite polyhedrally normed cone complex, without free faces, and with a convex bicombing. Then $X$ has finite combinatorial dimension and is CUB.

Let $X$ denote a locally finite polyhedrally normed complex without free faces. Then $X$ has at most one convex bicombing.
\ecor

\bp
Assume that $X$ is finite, and let $\sigma$ denote a convex bicombing. Then any face of $X$ is $\sigma$-convex. Since $X$ has no free face, for any distinct $x,y \in X$, the segment $\sigma(x,y)$ may be extended on each side. Since there are finitely many cells, we deduce that $\sigma$ is geodesically complete.

We can therefore apply Theorem~\ref{thm:descombes_lang_finite_dimension_unique_bicombing} and deduce that $X$ has finite combinatorial dimension and is CUB.

\mk

For the second part, assume that $\sigma$ is a convex bicombing on $X$. Then, according to the first part, $X$ is locally CUB. According to Theorem~\ref{thm:miesch_CUB_cartan_hadamard}, we deduce that $X$ is CUB, and $\sigma$ is the unique convex bicombing on $X$.
\ep

\subsection{Lattices and the CUB property}\

\mk

We will now explain how the results from Sections~\ref{subsec:horofunctions}~and~\ref{subsec:cone_complexes} are used in the proof of Theorem~\ref{thm:quotient_lattice_CUB} concerning the uniqueness of the convex bicombing on the diagonal quotient of the geometric realization of a lattice.

\mk

Let $L$ denote a bounded graded lattice of rank $n$, and let $X$ denote the geometric realization of $L$ with the standard $\ell^\infty$ metric.

Let $Y$ denote the geometric realization of the meet-semilattice $L \bs \{1\}$, with the standard polyhedral norm: it is called the \emph{diagonal quotient} of $X$.

Consider the affine version $M$ of $L$, with the diagonal isometric action of $\R$. Then the quotient metric space $\ov{M}=M / \R$ is called the \emph{diagonal quotient} of $M$.

\mk

We know the following. It is stated in \cite{haettel_helly_kpi1} for Euclidean buildings or Deligne complexes of Artin groups, but it relies in fact only on Theorem~\ref{thm:orthoscheme_lattice_injective} stating that the orthoscheme complex of a bounded graded lattice is injective.

\bthm[{\cite[Theorem~4.6]{haettel_helly_kpi1}}] \label{thm:orthoscheme_quotient_exists_convex_bicombing}
The metric space $\ov{M}$ has a convex bicombing $\sigma$. Moreover, the subspace $Y \subset \ov{M}$ is isometric and $\sigma$-convex. 
\ethm

We will improve that result to say that this convex bicombing is in fact unique.

\bthm \label{thm:quotient_lattice_CUB}
Let $L$ denote a bounded graded lattice of rank $n$, let $X$ denote the orthoscheme complex of $L$, and let $Y$ denote the diagonal quotient of $X$ endowed with the standard polyhedral metric. Then $Y$ is CUB.
\ethm

In order to prove this result, we first need to enlarge $L$ to ensure that the convex bicombing is geodesically complete.

\blem \label{lem:lattice_boolean_completion}
Let $L$ denote a bounded graded lattice of rank $n$. There exists a bounded graded lattice $L'$ of rank $n$ containing $L$ as a sublattice, such that each chain of $L'$ is contained in Boolean sublattice of $L'$.
\elem

\bp
For each maximal chain $c$ in $L$, let us consider a copy $B_{c}$ of the Boolean lattice ${\cal P}(\{1,\dots,n\})$, and consider the following poset $\theta(L)$:
$$\theta(L)= L \cup \left(\underset{\text{$c$ maximal chain in $L$}}{\bigcup} B_{c} /\sim \right),$$
where the chain $c \subset L$ is identified with the chain $(\emptyset,\{1\},\{1,2\},\dots,\{1,\dots,n\})$ in $B_{c}$.
Then $\theta(L)$ is a bounded graded lattice of rank $n$, which contains $L$ as a sublattice. Moreover, any chain in $L$ is contained in a Boolean sublattice of $\theta(L)$.

\mk

We can now define
$$L'=\bigcup_{j \in \N} \theta^j(L).$$
Then $L'$ is a bounded graded lattice of rank $n$, which contains $L$ as a sublattice. Moreover, any chain in $L'$ is contained in a Boolean sublattice of $L'$.
\ep

We will now prove a finiteness result for the combing lines in the affine version of a lattice.

\blem \label{lem:boolean_implies_geod_complete}
Let $L$ denote a bounded graded lattice of rank $n$, let $M$ denote the affine version of $L$, and let $\sigma$ denote the unique convex bicombing on $M$. Assume that any chain in $L$ is contained in a rank $n$ Boolean sublattice of $L$. Then $M$ is $\sigma$-geodesically complete. More precisely, for any faces $F,F'$ of $M$, there exists a finite subcomplex $Z$ of $M$ such that, for any distinct $x \in F$ and $y \in F'$, there exists an infinite $\sigma$-ray from $x$ through $y$ eventually contained in $Z$.
\elem

\bp
Let $x,y \in M$ distinct. Let $C$ denote a maximal cell containing a neighbourhood of $y$ in $\sigma(x,y)$, corresponding to a maximal chain $c$ of $L$. By assumption, there exists a Boolean sublattice $L_c$ of $L$ containing $c$. Let $M_c$ denote the subspace of $M$ corresponding to $L_c$: it is an isometric subspace of $M$, isometric to $\R^n$. As a consequence, $\sigma(x,y)$ can be extended infinitely inside $M_c$ as a $\sigma$-ray from $x$. Hence $M$ is $\sigma$-geodesically complete.
\ep

We are now able to gather all ingredients and prove Theorem~\ref{thm:quotient_lattice_CUB} stating that diagonal quotients of orthoscheme complexes have a unique convex bicombing. 

\bp[of Theorem~\ref{thm:quotient_lattice_CUB}]\

We will first prove that any finite subcomplex of $Y$ has finite combinatorial dimension. According to Lemma~\ref{lem:lattice_boolean_completion}, there exists a bounded graded lattice $L'$ of rank $n$, containing $L$ as a sublattice, such that each chain of $L'$ is contained in a Boolean sublattice of $L'$. Since the diagonal quotient of the orthoscheme complex of $L$ is isometrically embedded in the diagonal quotient of the orthoscheme complex of $L'$, we may assume that $L=L'$. More precisely, let us assume that any chain of $L$ is contained in a rank $n$ Boolean sublattice of $n$.

Let $M$ denote the affine version of $L$. Since $Y$ is isometrically embedded in $M$, it is sufficient to prove that any finite subcomplex of $M$ has finite combinatorial dimension. According to Theorem~\ref{thm:orthoscheme_quotient_exists_convex_bicombing}, there exists a convex bicombing $\sigma$ on $M$.

According to Lemma~\ref{lem:boolean_implies_geod_complete}, we know that the bicombing $\sigma$ on $M$ is geodesically complete. Moreover, $M$ satisfies the assumptions of Theorem~\ref{thm:geod_complete_implies_CUB}. As a consequence, we deduce that any finite subcomplex of $M$ has finite combinatorial dimension.

We deduce that any finite subcomplex of $Y$ has finite combinatorial dimension.

\mk

Remark that, since $M$ has finitely many shapes, it is locally conical. Now, according to Lemma~\ref{lem:finite_dimension_implies_unique_straight_geodesics}, we deduce that $Y$ is CUB and $\sigma$ is the unique convex bicombing on $Y$.
\ep

\section{Completion of posets} \label{sec:completion_of_posets}

We now show that we can relax the assumption that the (semi)lattice is graded, and replace it with the assumption that lengths of chains are uniformly bounded.

\bthm \label{thm:completion_semilattice}
Let $L$ denote a poset with minimum $0$, which is a meet-semilattice, with a bound on the length of chains. Assume that any $a,b,c \in L$ pairwise upperly bounded have a common upper bound. Then the orthoscheme complex of $L$, with the standard $\ell^\infty$ metric, is injective and CUB.
\ethm

\bp
Let $n \in \N$ be such that each chain of $L$ has length at most $n$. We will define a new poset $L_1$, containing $L$ as a subposet, which is graded. For each $x \in L$, let us define $r(x) \in \{0,1,\dots,n\}$ to be the length of a maximal chain from $0$ to $x$. Let us denote by $P$ the set of pairs of elements $(x,y) \in L^2$ such that $x<y$, $x$ is covered by $y$, and $r(y)-r(x) \geq 2$. For each $(x,y) \in P$, consider elements $u_{x,y,1},u_{x,y,2},\dots,u_{x,y,r(y)-r(x)-1}$. Let us now consider the poset
$$L_1 = L \cup \bigcup_{(x,y) \in P} \{u_{x,y,1},\dots,u_{x,y,r(y)-r(x)-1}\},$$
where the partial order on $L_1$ is generated by the order on $L$, and the relations $x<u_{x,y,1}<u_{x,y,2}<\dots<u_{x,y,r(y)-r(x)-1}<y$.
It is clear that $L_1$ is a graded poset of finite rank at most $n$, containing $L$ as a subposet.

\mk

Let $L_2$ denote a copy of $L_1$, and let $L'=L_1 \cup_L L_2$ denote the union of $L_1$ and $L_2$, where the two copies of the subposet $L$ are identified. Note that there is a natural involution $\theta$ on $L'$ exchanging $L_1$ and $L_2$ whose fixed point set is $L$.

\mk

The poset $L'$ has a minimum $0$, it is still a meet-semilattice, and is still such that any $a,b,c \in L'$ upperly bounded have a common upper bound in $L'$. Now $L'$ is also graded of rank at most $n$. According to Theorem~\ref{thm:orthoscheme_semilattice_injective}, the geometric realization of $L'$, with the orthoscheme $\ell^\infty$ metric, is injective and CUB.

\mk

Let $\sigma$ denote the unique convex bicombing on $|L'|$. By uniqueness, $\sigma$ is equivariant with respect to the involution $\theta$. In particular, the subspace $|L|$, which is the fixed-point set of $\theta$, is $\sigma$-convex. We deduce that $|L|$ admits a convex bicombing. Since $|L|$ is injective and finite-dimensional, according to Theorem~\ref{thm:descombes_lang_finite_dimension_unique_bicombing} we deduce that it is CUB.
\ep

\bthm \label{thm:completion_lattice_quotient}
Let $L$ denote a bounded lattice, with a bound on the length of chains. Let $X$ denote the orthoscheme complex of $L$, and let $Y$ denote the diagonal quotient of $X$ endowed with the standard polyhedral metric. Then $Y$ is CUB.
\ethm

\bp
The proof is quite similar to the proof of Theorem~\ref{thm:completion_semilattice}: let us write $L'=L_1 \cup_L L_2$, where $L_1, L_2$ are isomorphic bounded graded lattices containing $L$, and there is an involution $\theta$ on $L'$ exchanging $L_1$ and $L_2$ whose fixed point set is $L$.

\mk

The poset $L'$ is bounded and graded, so according to Theorem~\ref{thm:quotient_lattice_CUB}, the diagonal quotient $Y'$ of $|L'|$, endowed with the standard polyhedral metric, is CUB. Note that the involution $\theta$ induces an involution $\theta$ of $Y' = Y_1 \cup_Y Y_2$ whose fixed-point set is $Y$.

\mk

Let $\sigma'$ denote the unique convex bicombing on $Y'$. By uniqueness, $\sigma'$ is equivariant with respect to the involution $\theta$. In particular, the subspace $Y$, which is the fixed-point set of $\theta$, is $\sigma'$-convex. In particular, $Y$ admits a convex bicombing $\sigma$. As in the proof of Theorem~~\ref{thm:quotient_lattice_CUB}, we can prove that any finite subcomplex of $Y$ has finite combinatorial dimension, and according to Lemma~\ref{lem:finite_dimension_implies_unique_straight_geodesics} we conclude that $\sigma$ is the only convex bicombing on $Y$. As a consequence, $Y$ is CUB.
\ep

\section{Proof of the link conditions} \label{sec:proof_local_criteria}

We are now able to give a proof of the link conditions. Let us start with a technical lemma.

\blem \label{lem:nice_bowtie}
Let $L$ denote a poset with minimum $0$, with a bound on the length of chains. Assume that $L$ is not a meet-semilattice. Then there exists a bowtie $a,a' < b,b'$ in $L$ such that:
\bit
\item $a$ and $a'$ have a meet $m=a \wedge a'$.
\item The only element $c \in L$ such that $m \leq c \leq b,b'$ and $c$ is comparable to $a$ and $a'$ is $c=m$.
\eit
\elem

\bp
Let us call \emph{rank} of an element $x \in L$ the length of a maximal chain from $0$ to $x$. Since $L$ is not a meet-semilattice, according to Lemma~\ref{lem:characterizations_lattice_type_A}, $L$ contains a bowtie $a,a' < b,b'$. Among all bowties in $L$, assume that the ranks of $b,b'$ are minimal. Among all bowties $a,a' < b,b'$ such that the ranks of $b,b'$ are minimal, assume furthermore that the ranks of $a,a'$ are maximal.

\mk

In particular, since $a,a'$ have lower ranks than $b,b'$, they have a meet $m=a \wedge a'$.

\mk

Let us consider an element $c \in L$ such that $m \leq c \leq b,b'$ and $c$ is comparable to $a$ and $a'$. For a contradiction, assume that $m < c$. Since $m$ is the meet of $a,a'$, we deduce that $a<c$ or $a'<c$. For instance, assume that $a'<c$. We deduce that $a,c < b,b'$ is a bowtie in $L$, where the ranks of $a,c$ are greater than the ranks of $a,a'$: this is a contradiction. Hence $c=m$.
\ep

\mk

We are now able to give a proof of Theorem~\ref{thm:main_criterion_type_A_CUB}, for type-$A$ simplices.

\bp
Let $X$ denote a locally finite-dimensional flag simplicial complex with cyclically ordered simplices. Assume that, for every vertex $x \in X$, the set $(\St(x),\leq_x)$ is a meet-semilattice.

\mk

Fix a vertex $x \in X$, and consider the poset $L=(\St(x),\leq_x) \cup \{1\}$, where $1$ is the maximal element and $x$ is the minimal element. By assumption $L$ is a lattice. Furthermore, the length of chains of $L$ is bounded above by the maximal dimension of simplices of $X$ containing $x$. So we can apply Theorem~\ref{thm:completion_lattice_quotient} and deduce that the diagonal quotient $Y$ of $|L|$, with the standard polyhedral metric, is CUB. Since $X$ and $Y$ are locally isometric at $x$, we deduce that $X$ is locally CUB.

\mk

We now turn to the proof of the converse statement: if there exists a vertex $x \in X$ such that $(\St(x),\leq_x)$ is not a meet-semilattice, then $X$ is not locally CUB. According to Lemma~\ref{lem:nice_bowtie}, there exists a bowtie $a,a' <_x b,b'$ in $\St(x)$ such that:
\bit
\item $a$ and $a'$ have a meet $m=a \wedge_x a'$.
\item The only element $c \in \St(x)$ such that $m \leq_x c \leq_x b,b'$ and $c$ is comparable to $a$ and $a'$ is $c=m$.
\eit

\mk

For a contradiction, assume that $X$ has a conical bicombing $\sigma$. Note that the star of $b$ coincides with the $1$-ball centered at $b$, so the star of $b$ is $\sigma$-convex. Similarly, the stars of $b'$ and $m$ are $\sigma$-convex.

Since $a,a' <_x b,b'$ is a bowtie, $a$ and $a'$ are not comparable, then $a$, $a'$ do not lie in a common cell. Let us consider the two minimal simplices $C,C'$ of $\St(x)$ containing the beginning of $\sigma(a,a')$ starting from $a,a'$ respectively. By convexity, we know that $C$ and $C'$ are each adjacent to $b,b',m$. Let us consider $c$ the maximal element of $C \cap C'$. We know that $m \leq_x c \leq_x b,b'$, and $c$ is comparable to $a$ and $a'$. By assumption, we know that $c = m$.

Let $C_a$ denote the codimension $1$ face of $C$ opposite $a$: we have $d(a,C_a) = 1$, and $\sigma(a,a')$ exists $C$ through $C_a$. Similarly, let $C'_{a'}$ denote the codimension $1$ face of $C'$ opposite $a'$: we have $d(a',C'_{a'}) = 1$, and $\sigma(a,a')$ exists $C'$ through $C'_{a'}$. Since $d(a,a') \leq d(a,m)+d(m,a')=2$, we deduce that $\sigma(a,a')$ intersects $C_a \cap C'_{a'}$.

Let $q=\sigma(a,a',\f{1}{2}) \in C_a \cap C'_{a'}$: we have $q \leq_x m$ . We deduce that $\sigma(a,a') = [a,q] \cup [q,a']$. Now let $p \in [a,b]$ denote the midpoint of the edge $[a,b]$. We have $d(p,a)=\f{1}{2}$, $d(p,q)=1$ and $d(p,a') \leq 1$, which contradicts the conicality property.
\ep

We now turn to the proof of Theorem~\ref{thm:main_criterion_type_C_CUB}, for type $C$-simplices.

\bp

Let $X$ denote a locally finite-dimensional flag simplicial complex with ordered simplices, and assume that for every vertex $x \in X$, the poset $(\St(x),\leq_x)$ has no bowtie and, for any $a,b,c \in \St(x)$ which are pairwise upperly bounded (resp. lowerly bounded), they have a commun upper bound (resp. lower bound).

\mk

Endow $X$ with the standard $\ell^\infty$ metric, and fix a vertex $x \in X$. According to Proposition~\ref{pro:orthoscheme_complex_local_product}, a neighbourhood of $x$ in $X$ is isometric to a neighbourhood of $x$ in the $\ell^\infty$ product $|\St^+(x)| \times |\St^-(x)|$, where $\St^+(x)=\{y \in \St(x) \st x \leq_x y\}$ and $\St^-(x)=\{y \in \St(x) \st x \geq_x y\}$.

\mk

The poset $(\St^+(x),\leq_x)$ is a meet-semilattice with minimum $x$, whose length of chains is bounded above by the dimension of simplices of $X$ containing $x$. Moreover, any three pairwise upperly bounded elements have a common upper bound. According to Theorem~\ref{thm:completion_semilattice}, the geometric realization $|\St^+(x)|$ is injective and CUB.

Similarly the geometric realization $|\St^-(x)|$ is injective and CUB.

\mk

We deduce that the $\ell^\infty$ product $|\St^+(x)| \times |\St^-(x)|$ is injective and CUB, so in particular a neighbourhood of $x$ in $X$ is injective and CUB. Finally $X$ is locally injective and locally CUB.

\mk

We now turn to the proof of the converse statement: assume first that there exists a vertex $x \in X$ such that $(\St^+(x),\leq_x)$ has a bowtie. Then, by arguing as in the proof of Theorem~\ref{thm:main_criterion_type_A_CUB}, we conclude that $X$ has no convex bicombing.

\mk

Assume now that there exists a vertex $x \in X$ such that $(\St(x),\leq_x)$ has no bowtie, but there exist $a,b,c \in \St^+(x)$ which are pairwise upperly bounded, but have no common upper bound. We will prove that $X$ has no convex bicombing.

\mk

For a contradiction, assume that $X$ has a convex bicombing $\sigma$. Note that $\sigma(a,b)=[a,m] \cup [m,b]$, where $m$ denotes the midpoint of the edge $[a \wedge b,a \vee b]$. We also have $\sigma(a,c)=[a,m'] \cup [m',c]$, where $m'$ denotes the midpoint of the edge $[a \wedge c,a \vee c]$. Since $d(b,c) =1$ and $\sigma$ is convex, we deduce that $d(m,m') \leq \f{1}{2}$. Since $m,m'$ are midpoints of edges, this implies that $m,m'$ are in a cell of $\St(x)$. The maximal element $d$ of such a cell is such that $a \vee b,a \vee c \leq d$, so $d$ is an upper bound to $a,b,c$: this is a contradiction.
\ep

Finally, we turn to the proof of Theorem~\ref{thm:main_criterion_Garside_flag_CUB} for Garside flag complexes.

\bp
Let $(X,\varphi)$ denote a Garside flag complex, and endow $X$ with the standard $\ell^\infty$ metric. Note that $k$-simplices of $X$ come in columns, whose vertices come in chains of the form
$$\dots < \varphi^{-1}(x_k) < x_1 < x_2 < \dots < x_k < \varphi(x_1) < \varphi(x_2) < \dots.$$

\mk

In this sequence, any $k+1$ consecutive vertices form the vertices of a $k$-simplex of this column. In particular, one sees that there is an isometric action $(f_t)_{t \in \R}$ of $\R$ on $X$, such that $f_1=\varphi$.

As a consequence, given any point $p$ of $X$, up to translating $p$ using the action of $\R$ to a generic point, there exists a vertex $x \in X$ such that $p$ is contained in the interior of the geometric realization of the interval $[x,\varphi(x)]$. Since this interval is a homogeneous lattice, according to Theorem~\ref{thm:completion_semilattice}, we deduce that its geometric realization is injective and CUB. So $X$ is locally injective and locally CUB.

\mk

Now consider the quotient $Y=X/\varphi$, endowed with the standard polyhedral metric. Fix any vertex $y \in Y$, corresponding to the image of a vertex $x \in X$. Then a neighbourhood of $y$ in $Y$ is isometric to a neighbourhood of $y$ in the diagonal quotient of the lattice $[x,\varphi(x)]$. According to Theorem~\ref{thm:completion_lattice_quotient}, we deduce that a neighbourhood of $y$ in $Y$ is CUB. So $Y$ is locally CUB.
\ep

\section{Applications} \label{sec:applications}

We now show that the link conditions can be applied to numerous situations.

\subsection{Buildings}

We refer the reader to~\cite{abramenko_brown} and \cite{ronan} for references on buildings.

\mk

Consider a Euclidean building $X$ of type $\tilde{A}_n$: it is an $n$-dimensional simplicial complex, such that each vertex has a well-defined type in $\Z/(n+1)\Z$, and such that the vertices of each simplex have different types. Hence $X$ has cyclically ordered simplices. Moreover, for each vertex $x \in X$, the link $L$ of $x$ in $X$ is a spherical building of type $A_n$.

For instance, in case $X$ is the Euclidean building of $\SL(n+1,\Q_p)$, then $L$ is the spherical building of $\SL(n,\F_p)$. In other words, $L$ is the poset of non-trivial vector subspaces of $\F_p^n$, which is a lattice (up to adding $\{\{0\},\F_p^n\}$). 

More generally, $L$ is the poset of non-trivial subspaces of a projective geometry, and it is always a lattice (up to adding $\{0,1\}$). According to Theorem~\ref{thm:main_criterion_type_A_CUB} we deduce the following.

\bthm
Any Euclidean building $X$ of type $\tilde{A}_n$, with the standard polyhedral metric, is CUB.
\ethm

We can also consider a Euclidean building $X$ of extended type $\tilde{A}_n$: it is an $(n+1)$-dimensional simplicial complex, such that simplices have a well-defined total order. Similarly, we have the following.

\bthm
Any extended Euclidean building $X$ of type $\tilde{A}_n$, with the standard $\ell^\infty$ metric, is CUB and injective.
\ethm

\bp
Given an extended Euclidean building $X$ of type $\tilde{A}_n$, there is a well-defined canonical automorphism $\varphi$ of $X$, see~\cite{hirai_uniform_modular}. For instance, if $X$ is the Bruhat-Tits building of $\GL(n,\Q_p)$, then $\varphi$ is the homothety of $p$. For each vertex $x \in X$, the interval $[x,\varphi(x)]$ is the lattice of subspaces of a projective geometry of dimension $n$.

Then $(X,\varphi)$ is a Garside flag complex (see Definition~\ref{def:garside_flag_complex}): according to Theorem~\ref{thm:main_criterion_Garside_flag_CUB}, we deduce that $X$, with the standard $\ell^\infty$ metric, is CUB and injective.
\ep

Consider a Euclidean building $X$ of type $\tilde{B}_n$, $\tilde{C}_n$ or $\tilde{D}_n$. We will explain below how, up to subdividing simplices, the simplicial complex $X$ can be considered as a (non-thick) Euclidean building of type $\tilde{C}_n$ when we focus only on the simplicial structure.

\mk

Let us first consider the case of a classical spherical building. As explained in~\cite[Section~6.7]{abramenko_brown}, if $\K$ is a field of characteristic different from $2$, there are two natural spherical buildings associated with the group $\SO_{2n}(\K)$: the first one consists of the flag complex of nonzero totally isotropic subspaces of $\K^{2n}$, it has type $C_n$ and is not thick. The other one is the oriflamme complex of isotropic subspaces, it has type $D_n$ and is thick. The first building is just a simplicial subdivision of the second one.

\mk

More generally, the Euclidean Coxeter group of type $\tilde{B_n}$ or $\tilde{D_n}$ is a reflection subgroup of the Euclidean Coxeter group of type $\tilde{C_n}$. Indeed, the hyperplanes arrangements of types $\tilde{B_n}$ or $\tilde{D_n}$ are subarrangements of an arrangement of type $\tilde{C_n}$. For instance, as explained in~\cite[Section~6.7]{davis_coxeter} in the finite case, a model for the reflection hyperplanes in type $\tilde{D_n}$ is given by the hyperplanes $\{x_i \pm x_j =k\}_{1 \leq i<j \leq n, k \in \Z}$ in $\R^n$, while a model for the reflection hyperplanes in type $\tilde{C_n}$ is given by the hyperplanes $\{x_i \pm x_j =k\}_{1 \leq i<j \leq n, k \in \Z} \cup \{x_i =k\}_{1 \leq i \leq n, k \in \Z}$ in $\R^n$.

\mk

So if $X$ is a Euclidean building of type $\tilde{B}_n$ or $\tilde{D}_n$, one may subdivide each simplex into $2$ simplices (in type $\tilde{B}_n$) or $4$ simplices (in type $\tilde{D}_n$) so that each apartment is the Coxeter complex of type $\tilde{C_n}$, thus turning $X$ into a non-thick Euclidean building of type $\tilde{C_n}$.

\mk

Now let us consider a Euclidean building $X$ of type $\tilde{C_n}$: it is an $n$-dimensional simplicial complex, such that each vertex has a well-defined type in $\{0,1,\dots,n\}$, and the vertices of each simplex have different types. So $X$ has ordered simplices.

Moreover, for each vertex $x \in X$ of type $k$, the link $L$ of $x$ in $X$ is the join of two spherical buildings of types $B_k$ and $B_{n-k}$.

Then $L$ is the poset of non-trivial subspaces of a polar geometry: so $L \cup \{0\}$ is a meet-semilattice, and if three elements have a pairwise upper bound, they have a global upper bound. According to Theorem~\ref{thm:main_criterion_type_C_CUB}, we deduce the following.

\bthm
Any Euclidean building $X$ of type $\tilde{B}_n$, $\tilde{C}_n$ or $\tilde{D}_n$, with the standard $\ell^\infty$ metric, is CUB.
\ethm

\subsection{Simplices of groups}

We refer the reader to~\cite{bridson_haefliger} for the theory of complexes of groups. The link condition from Theorem~\ref{thm:main_criterion_type_A_CUB} can be used to prove that certain simplices of groups are developable, and that their development cover is CUB, as noticed in Theorem~\ref{thm:complex_groups_CUB_developable}.

\mk

Consider a simple $(n-1)$-dimensional simplex of groups $S$, with cyclically ordered vertices $(s_i)_{i \in \Z/n\Z}$. For each non-empty subset $I \subset \Z/n\Z$, let $G_I$ denote the group associated to the face $I$. When $I \subset J$, we will consider $G_J$ as a subgroup of $G_I$.

\bthm \label{thm:simplex_of_groups}
Assume that the simplex of groups satisfies the following.
\bit
\item Given any $i \in \Z/n\Z$ and $I,J \subset \Z/n\Z$ containing $i$, we have $G_I \cap G_J = G_{I \cup J}$.
\item Given any $i<j<k<\ell$ in $\Z/n\Z$, we have $G_{ik} \subset G_{ij}G_{i\ell}$.
\item Given any $i<j<k$ in $\Z/n\Z$, if $a,a' \in G_{ij}$ and $b,b' \in G_{ik}$ are such that $aba'b'=e$, we have $a,a' \in G_{ik}$ or $b,b' \in G_{ij}$ or there exists $i<j<\ell<k$ such that $a,a',b,b' \in G_{i\ell}(G_{ij} \cap G_{ik})$.
\eit
Then the simplex of groups is developable. Moreover its development, endowed with the standard polyhedral norm, is CUB.
\ethm

\bp
Since $S$ is a simple complex of groups, we know that $S$ is locally developable at each vertex. Fix a vertex $G_i$, and we will consider the local development $L_i$ at $G_i$: according to the first condition, it is the flag simplicial complex with vertex set $\{G_i\} \cup \bigcup_{j \neq i} G_i \bs G_{ij}$, where $aG_{ij}$ is adjacent to $bG_{ik}$ if and only if there exists $c \in G_i$ such that $cG_{ij}=aG_{ij}$ and $cG_{ik}=bG_{ik}$. Note that edges of $L_i$ have a well-defined orientation, where $aG_{ij} < bG_{ik}$ if $i<j<k$ (and $G_i$ is the minimal element of $L_i$).

\mk

We will first check that $L_i$ is the geometric realization of the corresponding poset. Assume that in $L_i$ we have two consecutive ordered edges between $aG_{ij}$, $bG_{ik}$ and $cG_{i\ell}$, where $a,b,c \in G_i$ and $i<j<k<\ell$. Without loss of generality, we may assume up to translation that $a=e$. Since $G_i,G_{ij},bG_{ik}$ form a simplex, we have $b \in G_{ij}$, so we may assume up to translation that $b=e$. Since $G_i,G_{ik},cG_{i\ell}$ form a simplex, we have $c \in G_{ik}$.

By assumption, we have $c = xy$, where $x \in G_{ij}$ and $y \in G_{i\ell}$. Hence $(G_i,G_{ij},cG_{i\ell})=x \cdot (G_i,G_{ij},G_{i\ell})$, so it is a simplex. Hence there is an ordered edge between $aG_{ij}$ and $cG_{i\ell}$ in the star of $G_i$. So $L_i$ is the geometric realization of the corresponding poset.

\mk

We will now check that the vertex set of $L_i$ is a meet-semilattice. According to Proposition~\ref{pro:lattice_balanced_bowtie}, we only need to check that the vertex set of $L_i$ has no balanced bowtie. Assume that there are $i<j<k$ and vertices such that $gG_{ij},gabG_{ij} < gaG_{ik},gaba'G_{ik}$, where $g \in G_i$, $a,a' \in G_{ij}$ and $b \in G_{ik}$. Up to translation, there exists $b' \in G_{ik}$ such that $aba'b'=e$ and the vertices are $G_{ij}=aba'b'G_{ij},abG_{ij} < aG_{ik},aba'G_{ik}$. 

\mk

By assumption, there are three possibilities. If $a,a' \in G_{ik}$, then $aG_{ik} = aba'G_{ik}$, so it is not a bowtie. If $b,b' \in G_{ij}$, then $G_{ij} = abG_{ij}$, so it is not a bowtie. Assume then that there exists $i<j<\ell<k$ such that $a,a',b,b' \in G_{i\ell}(G_{ij} \cap G_{ik})$. So the element $G_{i\ell}$ is such that $G_{ij}=aba'b'G_{ij},abG_{ij} \leq G_{i\ell} \leq aG_{ik},aba'G_{ik}$. This shows that there is no balanced bowtie in $L_i$.

\mk

According to Proposition~\ref{pro:lattice_balanced_bowtie}, we deduce that the vertex set of $L_i$ is a meet-semilattice. So according to Theorem~\ref{thm:main_criterion_type_A_CUB}, we see that the local development at each vertex is CUB. According to Theorem~\ref{thm:complex_groups_CUB_developable}, we deduce that the simplex of groups $S$ is developable. Moreover its development, endowed with the standard polyhedral norm, is CUB according to Theorem~\ref{thm:miesch_CUB_cartan_hadamard}.
\ep

Note that in the case of a triangle of groups, the criterion is the same as the one to ensure nonpositive curvature (with the equilateral triangle Euclidean norm), namely that the links of local developments have girth at least $6$.

However, already in the case of a $3$-simplex of groups, this is to our knowledge the first general combinatorial criterion ensuring developability. Here are very simple examples of such $3$-simplices of groups where the assumptions hold.

\bexe
For each $i \in \Z/4\Z$, let $G_i$ be the symmetric group $\frak{S}_4$. For each $j \in \{1,2,3\}$, define the image of $G_{i,i+j}$ in $G_i$ to be the stabilizer of $\{1,\dots,j\}$. For each $I \subset \{1,\dots,4\}$ containing $i$, define the image of $G_I$ in $G_i$ to be the intersection of all $G_{i,j}$, for $j \in I$. Then the conditions of Theorem~\ref{thm:simplex_of_groups} are easily seen to hold, and in fact the development of $S$ is the standard $\tilde{A}_3$ tiling of $\R^3$.
\eexe

\bexe
One can describe the Bruhat-Tits building of $\SL(n,\K((X)))$, for $n \geq 2$ and for any field $\K$, as follows. Consider the $(n-1)$-simplex of groups $T$ with vertices indexed by $\Z/n\Z$. For each $i \in \Z/n\Z$, let $t_i \in \GL(n,\K(X))$ denote the diagonal matrix with the first $i$ entries equal to $X$, and the last $n-i$ entries equal to $1$ (note that since we will make $t_i$ act by conjugation, this will be well-defined for $i \in \Z/n\Z$). For each $i \in \Z/n\Z$, let $G_i=t_i \SL(n,\K[X]) t_i^{-1} \subset \SL(n,\K(X))$. For each non-empty $I \subset \Z/n\Z$, let $G_I = \cap_{i \in I} G_i$. Note that, for each $i \in \Z/n\Z \bs \{0\}$, one can also describe the image of $G_{0,i}$ in $G_0$ as the stabilizer of the canonical $i$-plane $\K e_1 + \dots + \K e_i$ under the quotient action of $G_0=\SL(n,\K[X]) \ra \SL(n,\K)$ on $\K^n$. Then the conditions of Theorem~\ref{thm:simplex_of_groups} are easily seen to hold, and in fact the development of $T$ is precisely the Bruhat-Tits building of $\SL(n,\K((X)))$.
\eexe

\bexe
Let $A$ denote the Artin group of affine type $\tilde{A}_{n-1}$, with standard generators $(s_i)_{i \in \Z/n\Z}$. For each non-empty $I \subset \Z/n\Z$, let $G_I$ denote the standard parabolic subgroup generated by $\{s_i \st i \not\in I\}$. Then the corresponding $(n-1)$-simplex of groups $T$ satisfies the assumptions of Theorem~\ref{thm:simplex_of_groups}, and its development is the Artin complex of $A$ studied in Theorem~\ref{thm:artin_complex_An_CUB}.
\eexe

\subsection{Weak Garside groups}

We refer the reader to~\cite{garside}, \cite{bessis}, \cite{bessis2006garside} and \cite{haettel_huang_weakly_modular} for references concerning Garside and weak Garside groups. 

\mk

The classical definition of a Garside group starts with the definition of Garside monoid, as follows.

\bdf
A \emph{Garside monoid} (sometimes called quasi-Garside) is a pair $(M,\Delta)$, where $M$ is a monoid and
\bit
\item $M$ is left and right cancellative,
\item there exists $r:M \bs \{1\} \ra \N \bs \{0\}$ such that $r(fg) \geq r(f)+r(g)$,
\item any two elements of $M$ have left and right least common multiples and greatest common divisors,
\item $\Delta$ is a \emph{Garside} element of $M$, i.e. the family $S$ of left and right divisors of $\Delta$ coincide and generate $M$.
\eit 
If one further require that the set $S$ of simple elements is finite, $(M,\Delta)$ is called of finite type.
\edf

Now we can define a Garside group.

\bdf
A group $G$ is called a \emph{Garside group} if there exists a Garside monoid $(M,\Delta)$ such that $G$ is the group of left fractions of $M$.
\edf

\mk

Examples of Garside groups of finite type include finite rank free abelian groups, braid groups, and more generally spherical-type Artin groups. Tree products of cyclic groups also form a nice family of examples, see~\cite{picantin_amalgams_Garside}. The standard Garside structure for a spherical-type Artin group is associated to the standard Artin monoid.  Note that, according to Bessis (\cite{bessis}), spherical-type Artin groups actually admit another Garside structure, called the dual Garside structure. See notably~\cite{paolini_dual_approach},\cite{paolini_salvetti}, \cite{haettel_bourbaki_paolini_salvetti} and \cite{delucchi_paolini_salvetti_rank_3} for the importance of the dual approach for the study of Artin groups.

\mk

Examples of Garside groups (with an infinite set of simple elements) include finite rank free groups (see~\cite{bessis_free}) and Euclidean Artin groups of type $\tilde{A}_n$, $\tilde{C}_n$ and $\tilde{G}_2$ (see~\cite{brady_mccammond_factoring}, \cite{mccammond_sulway} and \cite{mccammond_failure_lattice_property}). Note that the dual structure of Euclidean-type Artin groups is not Garside; nevertheless, McCammond and Sulway have defined a way to "complete" these groups to obtain crystallographic braid groups (\cite{mccammond_sulway}), which are fundamental in the proof of the $K(\pi,1)$ conjecture for Euclidean-type Artin groups by Paolini and Salvetti (\cite{paolini_salvetti}).

\mk

Another closely related notion is that of weak Garside group. Rather than giving the definition, we prefer stating the following geometric characterization of Garside and weak Garside groups.

\bthm[{\cite[Theorem~4.7]{haettel_huang_weakly_modular}}] \label{thm:dictionary_garside_group}
A group $G$ is a Garside group (resp. weak Garside group) if and only if there exists a Garside flag complex $(X,\varphi)$ such that $G$ can be realized as a group of automorphisms of $X$ commuting with $\varphi$, acting freely and transitively (resp. freely) on vertices of $X$.

Moreover, the group $G$ is a (weak) Garside group of finite type if and only if $X$ can be chosen such that the action of $G$ is cocompact.\ethm

For instance, finite index subgroups of Garside groups are weak Garside groups. Other examples include the fundamental groups of complements of real simplicial arrangements of hyperplanes (\cite{deligne}), all braid groups of complex reflection groups (\cite{bessis_noncrossing_partitions_complex},\cite{bessis_finite_complex_Kpi1}, \cite{corran_new_garside_braid}) except possibly the exceptional complex braid group of type $G_{31}$, and some extensions of Artin groups of type $B_n$ (\cite{crisp_paris_representations_braid}).

\mk

As a consequence of Theorem~\ref{thm:main_criterion_Garside_flag_CUB}, we have the following.

\bcor
Let $(G,\Delta)$ denote a weak Garside group of finite type, and let $k \geq 1$ be such that $\Delta^k \in Z(G)$. Then $G$ and $G/\<\Delta^k\>$ both act properly and cocompactly by isometries on a CUB space.
\ecor

\bp
If we denote by $(X,\varphi)$ the Garside flag complex associated to $G$, we know that $G$ acts geometrically on $X$, and $G/\<\Delta^k\>$ acts geometrically on $X/ \varphi$.
\ep

Let $A$ be a spherical-type Artin group and $\Delta$ be a Garside element of $A$. If we denote by $(X,\varphi)$ the associated Garside flag complex, then $X/\varphi$ has also been described by Bestvina (\cite{bestvina_artin}) as the normal form complex of $A$, which exhibits some form of combinatorial nonpositive curvature. We strengthen this claim by remarking that, when we endow Bestvina's complex with the standard polyhedral metric, it is CUB.

\subsection{Artin complexes for some Euclidean-type Artin groups}

We refer the reader to \cite{paris_kpi1}, \cite{paris_kpi1}, \cite{godelle_paris}, \cite{charney_davis}, \cite{charney_davis_kpi1}, \cite{charney_problems}, \cite{mccammond_mysterious} for references on Artin groups (also called Artin-Tits groups).

\mk

Let $A$ denote any Artin group. The Artin complex is the flag simplicial complex $X$ whose vertex set consists of left cosets of maximal proper standard parabolic subgroups of $A$, with an edge between $gP$ and $g'P'$ if and only if $gP \cap g'P' \neq \emptyset$, see~\cite{cumplido_martin_vaskou_parabolic} and \cite[Remark~(i), p.~606]{charney_davis_kpi1}. Note that, from the presentation of $A$, the complex $X$ is simply-connected.

We will be interested in the case where $A$ is of Euclidean-type $\tilde{A}_n$ or $\tilde{C}_n$.

\mk

In type $\tilde{A}_n$, vertices of $X$ have a well-defined type in $\Z/(n+1)\Z$, so simplices of $X$ have a well-defined cyclic order, and the link $L$ of any vertex is isomorphic to the Artin complex of the Artin group of type $A_n$. Bessis (\cite{bessis_free}), and independently Crisp and McCammond (unpublished), proved that $L \cup \{0,1\}$ is isomorphic to the lattice of cut-curves (see~\cite{haettel_helly_kpi1} for a proof, following Crisp and McCammond). According to Theorem~\ref{thm:main_criterion_type_A_CUB}, we deduce the following.

\bthm \label{thm:artin_complex_An_CUB}
The Artin complex of type $\tilde{A_n}$, with the standard polyhedral metric, is CUB.
\ethm

\mk

In type $\tilde{C}_n$, vertices of $X$ have a well-defined type in $\{0,\dots,n\}$, so simplices of $X$ have a well-defined total order. Moreover, the link $L_k$ of any vertex of type $k$ is the join of the Artin complexes of types $B_k$ and $B_{n-k}$. According to~\cite[Proposition~6.3]{haettel_helly_kpi1}, the Artin complex of type $B_k$, with $0$, is a meet-semilattice such that any $a,b,c$ which are pairwise upperly bounded have an upper bound. According to Theorem~\ref{thm:main_criterion_type_C_CUB}, we deduce the following.

\bthm
The Artin complex of type $\tilde{C_n}$, with the standard $\ell^\infty$ metric, is CUB.
\ethm

\subsection{Some arc complexes}

The question of finding nonpositive curvature metrics on some curve complexes or arc complexes on surfaces is quite intriguing and difficult, raised notably by Masur and Minsky (\cite{masur_minsky}). For instance, Webb proved that many such complexes do not admit CAT(0) metrics (\cite{webb_non_CAT0}).

\mk

The following is nothing more than a topological description of the Artin complex of the Artin group of type $\widetilde A_{n-1}$ (see~\cite[Proposition~5.8]{haettel_huang_weakly_modular}). Let $\Sigma$ denote a $2$-sphere with $n+2$ punctures $\{N,S,p_1,\dots,p_n\}$ with two distinguished punctures $N,S$ which could be thought of as the North Pole and the South Pole of $\Sigma$. The punctures $p_1,\dots,p_n$ may be thought as cyclically ordered on the equator of the $2$-sphere.

\mk

\begin{figure}[H]
\centering
\includegraphics[width=5cm]{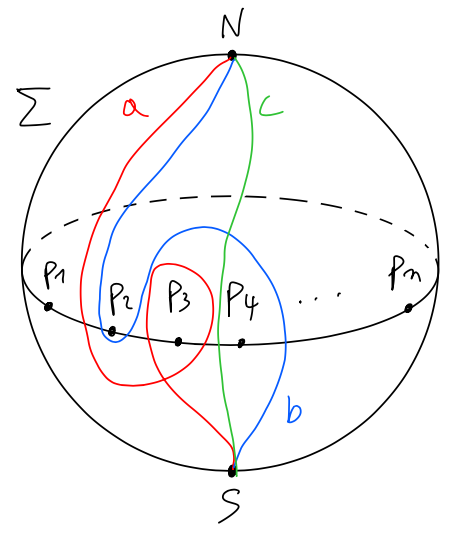}
\caption{Arcs on the punctured sphere $\Sigma$: $a$ is adjacent to $b$ and $c$}
\label{fig:arc_complex}
\end{figure}

Let ${\cal A}(\Sigma)$ denote the following simplicial complex. Its vertex set consists of isotopy classes of arcs in $\Sigma$ from $N$ to $S$. Two vertices are adjacent if they can be realized disjointly. Then ${\cal A}(\Sigma)$ is the associated flag simplicial complex, see Figure~\ref{fig:arc_complex}. According to~\cite[Lemma~2.5]{wahl}, this arc complex ${\cal A}(\Sigma)$ is contractible.

\mk

There is a total cyclic order on vertices of simplices of ${\cal A}(\Sigma)$, when one fixes an orientation on $\Sigma$: let $\sigma$ denote a simplex of ${\cal A}(\Sigma)$. Then any two distinct $a,b \in \sigma$ are such that $\Sigma \bs \{a,b\}$ has two connected components, associated to a fixed orientation of $\Sigma$. So given any three $a,b,c \in \sigma$, we say that $a<b<c$ if $a$ is on the left of $b$, and $c$ is on the right of $b$.

\bthm
The complex ${\cal A}(\Sigma)$, with the standard polyhedral norm, is CUB.
\ethm

\bp
Fix an arc $a \in {\cal A}(\Sigma)$, and let us consider the surface $D=\Sigma \bs a$: it is homeomorphic to a disk with $n$ punctures, and with two marked points $N,S$ on the boundary $\partial D$. The link $L$ of $a$ in ${\cal A}(\Sigma)$ is isomorphic to the complex of arcs in $D$ from $N$ to $S$. Let us call $a_W,a_E$ the two connected components of $\partial D \bs \{N,S\}$, where $E$ stands for East and $W$ for West.

\mk

\begin{figure}[H]
\centering\includegraphics[width=5cm]{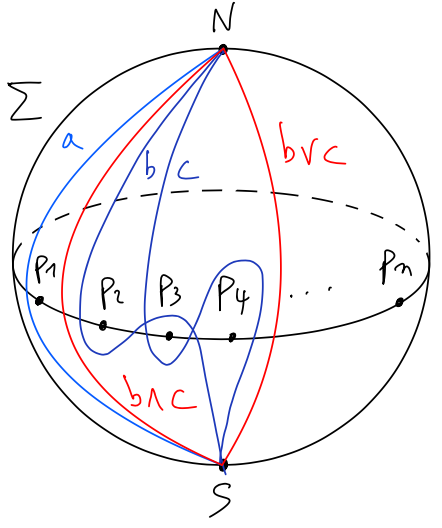}
\caption{The lattice property: the meet $b \wedge c$ and the join $b \vee c$ of the two arcs $b,c$ in the interval $[a_W,a_E]$.}
\label{fig:lattice}
\end{figure}

The induced order on $L$ is the following: if $b,c$ are disjoints arcs in $D$ from $N$ to $S$, then $b \leq_a c$ if $b$ is on the west of $c$. The poset $L$ can be completed with a minimum element $a_W$, and a maximum element $a_E$. Now $L \cup \{a_W,a_E\}$ is isomorphic to the poset of cut-curves. According to~\cite{bessis_free} (and~\cite{haettel_helly_kpi1} for an account of the unpublished proof due to Crisp and McCammond), this poset is a lattice. See also~\cite{haettel_huang_weakly_modular}, or Figure~\ref{fig:lattice} to see how the lattice property works.

\mk

We can therefore apply Theorem~\ref{thm:main_criterion_type_A_CUB} and deduce that the complex ${\cal A}(\Sigma)$, with the standard polyhedral norm, is CUB.
\ep

\subsection{A complex of homologous multicurves}

Apart from the complex of arcs on a punctured sphere, there is another natural complex for which the lattice property is straightforward. Note that the complex, as well the tools used for its study, are very similar to the Kakimizu complex studied in~\cite{przytycki_schultens_kakimizu}, see also Section~\ref{subsec:kakimizu}.

\mk

Fix a closed surface $S$ of genus $g \geq 1$ with $p \geq 0$ punctures. Fix the homology class $[a] \in H_1(S,\Z)$ of an oriented nonseparating simple closed curve. Let $\tilde S$ denote the associated infinite cyclic cover. 

Let us remind that a multicurve on $S$ is a finite union of disjoint, non-homotopic simple closed curves. If $b$ is an oriented multicurve on $S$ homologous to $a$, the lifts of the complement $S \bs b$ to $\tilde S$ are separated by lifts of $b$ denoted $(\tilde b_n)_{n \in \Z}$. Also note that each $\tilde b_n$ is a boundary component of an unbounded subsurface $B(\tilde b_n)$ of $\tilde S$ which is "below" $\tilde{b_n}$, i.e. which contains all lifts $\tilde b_m$ with $m < n$. If $b,b'$ are two such oriented multicurves on $S$ homologous to $a$, we say that a lift $\tilde b_n$ of $b$ is \emph{below} of a lift $\tilde b'_m$ of $b'$ if $B(\tilde b_n)$ is contained in $B(\tilde b'_m)$, up to homotopy. We denote $\tilde b_n \prec \tilde b'_m$ in this case.

\mk

Let ${\cal C}_a(S)$ denote the flag simplicial complex defined as follows:
\bit
\item Vertices of ${\cal C}_a(S)$ are homotopy classes of oriented multicurves on $S$ homologous to $a$.
\item There is an edge between $b,b' \in {\cal C}_a(S)$ if there are lifts $\tilde b_n$, $\tilde b'_m$ of $b,b'$ respectively, such that $\tilde b_n \prec \tilde b'_m \prec \tilde b_{n+1}$.
\eit

\begin{figure}[H]
\centering
\includegraphics[width=10cm]{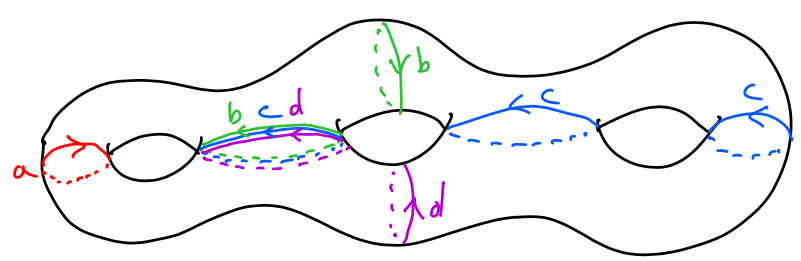}
\caption{Four pairwise disjoint multicurves $a,b,c,d$.}
\label{fig:simplex_multicurves}
\end{figure}

Here is an example of a closed surface $S$ of genus $3$, with a nonseparating simple closed curve $a$. There are also three disjoint multicurves $b,c,d$ represented, all of which are homologous to $a$, see Figure~\ref{fig:simplex_multicurves}.

\begin{figure}[H]
\centering
\includegraphics[width=10cm]{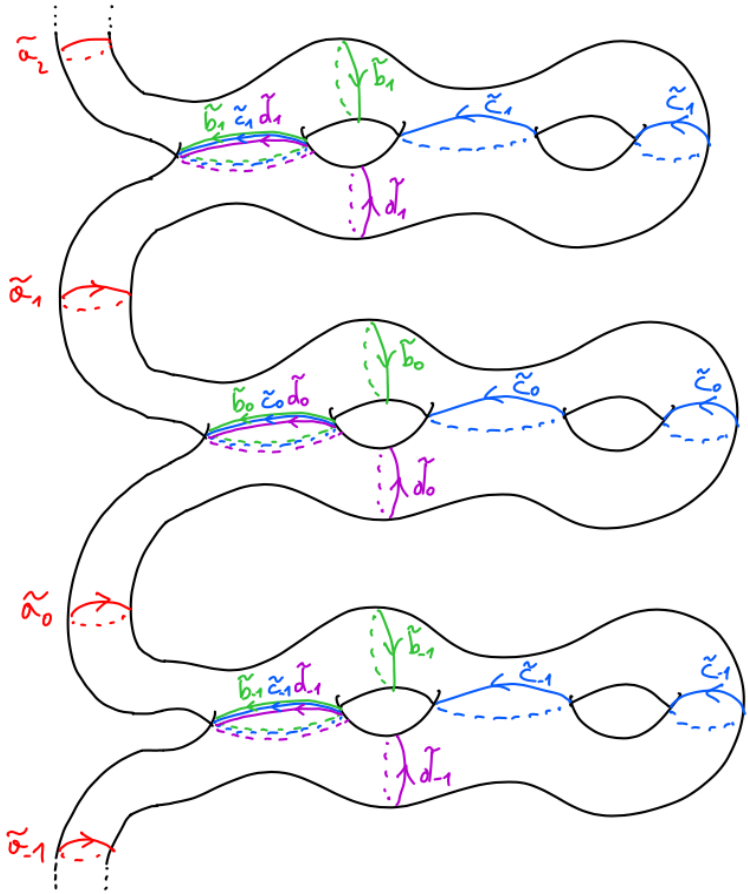}
\caption{The lifts of the multicurves $a,b,c,d$ in the cyclic cover.}
\label{fig:simplex_cyclic_cover}
\end{figure}

When we consider the cyclic cover $\tilde S$ of $S$ associated to $[a]$, we obtain a family of lifts of each of the multicurves $a,b,c,d$. In this example, the lifts satisfy $\tilde a_0 \prec \tilde d_0 \prec \tilde c_0 \prec \tilde b_0 \prec \tilde a_1$, so the four multicurves $a,b,c,d$ form a simplex of ${\cal C}_a(S)$, see Figure~\ref{fig:simplex_cyclic_cover}.

Hatcher and Margalit defined a very similar complex, and proved that it is contractible (\cite[Proposition~7]{hatcher_margalit_torelli}). We will prove that the complex ${\cal C}_a(S)$ is contractible, following ideas from~\cite{przytycki_schultens_kakimizu}.

\bpro
The complex ${\cal C}_a(S)$ is contractible.
\epro

\bp
Let us introduce the following distance on vertices of ${\cal C}_a(S)$. If $b,b' \in {\cal C}_a(S)$, the distance $D(b,b')$ is the minimal number $R \in \N$ such that there exist lifts $\tilde b_n$, $\tilde b'_m$ of $b,b'$ respectively, such that $\tilde b_n \prec \tilde b'_m \prec \tilde b_{n+R}$.

\mk

Note that $D$ is bounded above by the graph distance. In particular, if $D(a,b)=1$, then $b$ is adjacent to $a$. We will prove that there is a retraction from $B_D(a,R+1)$ to $B_D(a,R)$ for $R \in \N$.

\mk

Fix $b \in B_D(a,R+1)$. Consider the lifts $(\tilde a_n)_{n \in \Z}$, and consider the unique lift $\tilde b$ of $b$ such that $\tilde a_0 \prec \tilde b \prec a_{R+1}$. Let us consider the subsurface $T=B(\tilde b) \cap B(\tilde a_R)$ of $\tilde S$ up to homotopy: to be more precise, we may fix a hyperbolic structure on $S$, its lift to $\tilde S$, and only consider geodesic representatives of multicurves. The boundary $\tilde c_0=\partial T$ is a multicurve in $\tilde S$, whose projection $c$ to $S$ is a multicurve homologous to $a$. Moreover, we know that $\tilde c_0 \prec \tilde b \prec a_{R+1} \prec \tilde c_1$, so $c$ is adjacent to $b$. Moreover, since $\tilde a_0 \prec \tilde c_0 \prec a_R$, we know that $D(a,c) \leq R$.

\mk

So we have defined a simplicial retraction from $B_D(a,R+1)$ to $B_D(a,R)$ for $R \in \N$. As a consequence, the complex ${\cal C}_a(S)$ is contractible.
\ep

Note that if $\sigma$ is a simplex of ${\cal C}_a(S)$, we know that we can find lifts $\tilde{b^1_{n_1}},\dots,\tilde{b^p_{n_p}}$ of $\sigma$ to $\tilde S$ satisfying $\tilde{b^1_{n_1}} \prec \dots \prec \tilde{b^p_{n_p}} \prec \tilde{b^1_{n_1+1}}$. Furthermore, the corresponding cyclic order $b^1 < b^2 < \dots < b^p < b^1$ on vertices of $\sigma$ is well-defined. In the example of Figure~\ref{fig:simplex_multicurves}, the corresponding cyclic order is $a<b<c<d<a$. For this complex, we are able to prove the following.

\bthm
The complex ${\cal C}_a(S)$, endowed with the standard polyhedral metric, is CUB. In particular, it is contractible.
\ethm

\bp
Fix a multicurve $b \in {\cal C}_a(S)$, with a lift $\tilde b_0$ in $\tilde S$. Let us consider the set $E$ of lifts $\tilde c$ of elements $c \in {\cal C}_a(S)$ such that $\tilde b_0 \prec \tilde c \prec \tilde b_1$. The set $E$ may be endowed with the partial order $\prec$. We want to prove that the star of $b$, with the induced order, is a meet-semilattice. This poset is isomorphic to $(E \bs \{\tilde b_1\},\prec)$. So it is equivalent to prove that $(E,\prec)$ is a meet-semilattice.

\mk

Consider $\tilde c,\tilde c' \in E$. Let us consider the subsurface $T=B(\tilde c) \cap B(\tilde c')$ of $\tilde S$ up to homotopy. Its boundary $\tilde d=\partial T$ is the lift of a muticurve $d \in {\cal C}_a(S)$, such that $\tilde b_0 \prec \tilde d \prec \tilde c, \tilde c' \prec \tilde b_1$. In particular, $\tilde d \in E$. Moreover, by construction, $\tilde d$ is the maximal such element of $E$: this means that $\tilde d$ is the meet of $\tilde c$ and $\tilde c'$. So $E$ is a meet-semilattice.

\mk

\begin{figure}[H]
\centering
\includegraphics[width=10cm]{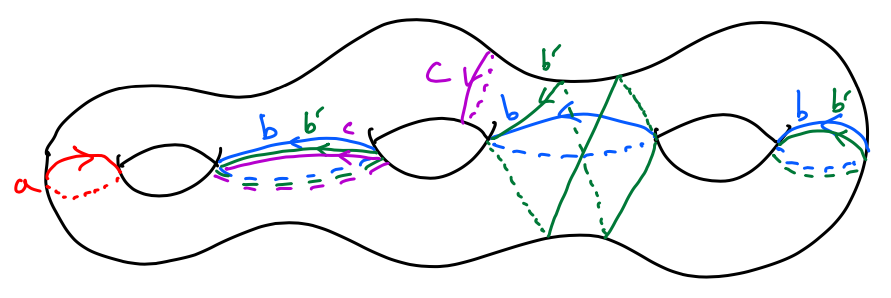}
\caption{The meet $c$ of the multicurves $b$ and $b'$.}
\label{fig:chirurgy}
\end{figure}

The meet of two multicurves in the star of $a$ is depicted in Figure~\ref{fig:chirurgy}. For simplicity, we represented the multicurves in the surface $S$ and not in its cyclic cover.

\mk

According to Theorem~\ref{thm:main_criterion_type_A_CUB}, we deduce that the complex ${\cal C}_a(S)$, endowed with the standard polyhedral metric, is locally CUB. Since the complex ${\cal C}_a(S)$ is simply-connected, we deduce by Theorem~\ref{thm:miesch_CUB_cartan_hadamard} that it is CUB.
\ep

Note that the stabilizer of the homology class of $a$ in the mapping class group of $S$ acts on ${\cal C}_a(S)$ by isometries. In particular, the Torelli group of $S$ acts on ${\cal C}_a(S)$ by isometries, for each nonseparating simple closed curve $a$ on $S$.

\subsection{The Kakimizu complex} \label{subsec:kakimizu}

Let $L$ denote a knot in $\SS^3$, and let $E=E(L)$ denote the exterior of a tubular neighbourhood of $L$. A \emph{spanning surface} is a surface properly embedded in $E$, which is contained in some Seifert surface for $L$. Let $MS(L)$ denote the flag simplicial complex whose vertices are isotopy classes of minimal genus spanning surfaces, with an edge between two surfaces if they have disjoint representatives.

\mk

This complex can also be defined for the complement of a link in $\SS^3$, or more generally for a compact, connected, orientable irreducible $3$-manifold with boundary, see~\cite{przytycki_schultens_kakimizu}.

\mk

This complex is called the \emph{Kakimizu complex}, and it has been defined in~\cite{kakimizu}. Scharlemann and Thompson proved that this complex is connected (\cite{scharlemann_thompson_kakimizu}), and Kakimizu gave another proof when $L$ is a link (\cite{kakimizu}). Schultens proved that $MS(L)$ is simply-connected (\cite{schultens_kakimizu}), and Przytycki and Schultens proved that it is in fact contractible (\cite{przytycki_schultens_kakimizu}).

\mk

Note that if we fix an orientation of $L$, there is a well-defined total cyclic ordering on vertices of simplices of $MS(L)$, given informally by the "angle" at which surfaces intersect the tubular neighbourhood $E$ of $L$. More precisely, consider the infinite cyclic covering $\tilde{E}$ of $E$ corresponding to the meridian of $L$. If $S,R,R'$ are pairwise disjoint spanning surfaces, then $R,R'$ have disjoint lifts $\tilde{R},\tilde{R'}$ in $\tilde{E}$ contained between two consecutives lifts $\tilde{S}_0,\tilde{S}_1$ of $S$. One then say that $S<R<R'$ if $\tilde{R}$ separates $\tilde{S}_0$ and $\tilde{R'}$.

\mk

Using this order, one may endow the Kakimizu complex with the standard polyhedral metric.

\bthm
The Kakimizu complex $MS(L)$, endowed with the standard polyhedral metric, is CUB.
\ethm

As a consequence, one obtains another proof that the Kakimizu is contractible (assuming that it is simply-connected). Applying Proposition~\ref{pro:action_CUB_fixed_point}, one also deduces another proof of~\cite[Corollary~1.3]{przytycki_schultens_kakimizu} that for any finite subgroup $G$ of the mapping class group of $E$, there exists a union of pairwise disjoint minimal genus spanning surfaces which is $G$-invariant up to isotopy. 

\bp
We will apply Theorem~\ref{thm:main_criterion_type_A_CUB}. Since $MS(L)$ is simply-connected according to~\cite{schultens_kakimizu}, it is sufficient to check that the star of a vertex is a meet-semilattice.

\mk

Fix three spanning surfaces $S,R,R' \in MS(L)$, with $R,R'$ disjoint from $S$, and consider the infinite cyclic covering $\tilde{E}$ of $E$ corresponding to the meridian of $L$. Let $\tau$ denote the covering automorphism of $\tilde{E}$. Fix lifts $\tilde{S},\tilde{R},\tilde{R'}$ of $S,R,R'$ to $\tilde{E}$ such that $\tilde{R},\tilde{R'}$ are disjoint from $\tilde{S}$, but not from $\tau^{-1}(\tilde{S})$. Let $\tilde{P}$ denote the surface "below $\tilde{R} \cup \tilde{R'}$", i.e. obtained from $\tilde{R}$ and $\tilde{R'}$ by a cut-and-paste operation as in~\cite[Section~4]{przytycki_schultens_kakimizu}. Let $P$ denote the image of $\tilde{P}$ in $E$. Then it is clear that we have $S \leq P \leq R,R'$. Furthermore, by construction $P$ is maximal, hence $P$ is the meet of $R$ and $R'$. So the star of $S$ is a meet-semilattice.
\ep

\section*{Conflict of interest}

The author states that there is no conflict of interest.

\section*{Data availability}

Data sharing not applicable to this article as no datasets were generated or analysed during the current study.

\bibliographystyle{alpha}
\bibliography{bibli}

\def\polhk#1{\setbox0=\hbox{#1}{\ooalign{\hidewidth
  \lower1.5ex\hbox{`}\hidewidth\crcr\unhbox0}}}
\begin{thebibliography}{HSWW17}

\bibitem[AB08]{abramenko_brown}
Peter Abramenko and Kenneth~S. Brown.
\newblock {\em {Buildings. Theory and Applications.\!}}
\newblock Grad. Text. Math. Springer, 2008.

\bibitem[Ago13]{agol}
Ian Agol.
\newblock The virtual {H}aken conjecture.
\newblock {\em Doc. Math.}, 18:1045--1087, 2013.
\newblock With an appendix by Agol, Daniel Groves, and Jason Manning.

\bibitem[Bal95]{ballmann_lectures_NPC}
Werner Ballmann.
\newblock {\em Lectures on spaces of nonpositive curvature}, volume~25 of {\em
  DMV Seminar}.
\newblock Birkh\"{a}user Verlag, Basel, 1995.
\newblock With an appendix by Misha Brin.

\bibitem[Bas18]{basso_fixed_point}
Giuliano Basso.
\newblock Fixed point theorems for metric spaces with a conical geodesic
  bicombing.
\newblock {\em Ergodic Theory Dynam. Systems}, 38(5):1642--1657, 2018.

\bibitem[Bas24]{basso_bicombings}
Giuliano Basso.
\newblock Extending and improving conical bicombings.
\newblock {\em Enseign. Math.}, 70(1-2):165--196, 2024.

\bibitem[BC06]{bessis_noncrossing_partitions_complex}
David Bessis and Ruth Corran.
\newblock Non-crossing partitions of type {$(e,e,r)$}.
\newblock {\em Adv. Math.}, 202(1):1--49, 2006.

\bibitem[Bes96]{bestvina_zstructure}
Mladen Bestvina.
\newblock Local homology properties of boundaries of groups.
\newblock {\em Michigan Math. J.}, 43(1):123--139, 1996.

\bibitem[Bes99]{bestvina_artin}
Mladen Bestvina.
\newblock Non-positively curved aspects of {A}rtin groups of finite type.
\newblock {\em Geom. Topol.}, 3:269--302, 1999.

\bibitem[Bes03]{bessis}
David Bessis.
\newblock The dual braid monoid.
\newblock {\em Ann. Sci. \'Ecole Norm. Sup. (4)}, 36(5):647--683, 2003.

\bibitem[Bes06a]{bessis_free}
David Bessis.
\newblock A dual braid monoid for the free group.
\newblock {\em J. Algebra}, 302(1):55--69, 2006.

\bibitem[Bes06b]{bessis2006garside}
David Bessis.
\newblock Garside categories, periodic loops and cyclic sets.
\newblock {\em arXiv preprint math/0610778}, 2006.

\bibitem[Bes15]{bessis_finite_complex_Kpi1}
David Bessis.
\newblock Finite complex reflection arrangements are {$K(\pi,1)$}.
\newblock {\em Ann. of Math. (2)}, 181(3):809--904, 2015.

\bibitem[BEW13]{bestvina_eskin_wortman}
Mladen Bestvina, Alex Eskin, and Kevin Wortman.
\newblock Filling boundaries of coarse manifolds in semisimple and solvable
  arithmetic groups.
\newblock {\em J. Eur. Math. Soc. (JEMS)}, 15(6):2165--2195, 2013.

\bibitem[BGS85]{bgs}
Werner Ballmann, Mikhael Gromov, and Viktor Schroeder.
\newblock {\em {Manifolds of nonpositive curvature\!}}
\newblock {Progr.~Math. {\bf 61}}. {Birkh\"auser}, 1985.

\bibitem[BH99]{bridson_haefliger}
Martin~R. Bridson and Andr\'e Haefliger.
\newblock {\em {Metric \,spaces \,of \,non-positive \,curva\-ture}}, volume
  {319\!} of {\em {Grund.~math.~Wiss.}}
\newblock {Springer}, 1999.

\bibitem[BJ07]{borel_ji}
Armand Borel and Lizhen Ji.
\newblock {Compactifications of symmetric spaces}.
\newblock {\em {J.~Diff.~Geom.}}, {75}:1--56, {2007}.

\bibitem[BM10]{brady_mccammond}
Tom Brady and Jon McCammond.
\newblock Braids, posets and orthoschemes.
\newblock {\em Algebr. Geom. Topol.}, 10(4):2277--2314, 2010.

\bibitem[BM15]{brady_mccammond_factoring}
Noel Brady and Jon McCammond.
\newblock Factoring {E}uclidean isometries.
\newblock {\em Internat. J. Algebra Comput.}, 25(1-2):325--347, 2015.

\bibitem[Bow95]{bowditch}
Brian~H. Bowditch.
\newblock {Notes on locally {${\rm CAT}(1)$} spaces}.
\newblock In {\em Geometric group theory ({C}olumbus, {OH}, 1992)}, volume~3 of
  {\em Ohio State Univ. Math. Res. Inst. Publ.}, pages 1--48. 1995.

\bibitem[Bri06]{brill_horofunction_polyhedral}
Bj{\"o}rn Brill.
\newblock {Eine Familie von Kompaktifizierungen affine Geb{\"a}ude}.
\newblock {\em PhD thesis, Frankfurt}, 2006.

\bibitem[Bus48]{busemann_npc}
Herbert Busemann.
\newblock Spaces with non-positive curvature.
\newblock {\em Acta Math.}, 80:259--310, 1948.

\bibitem[BW02]{brady_watt}
Thomas Brady and Colum Watt.
\newblock {{$K(\pi,1)$}'s for {A}rtin groups of finite type}.
\newblock In {\em Proceedings of the {C}onference on {G}eometric and
  {C}ombinatorial {G}roup {T}heory, {P}art {I} ({H}aifa, 2000)}, volume~94,
  pages 225--250, 2002.

\bibitem[CCHO21]{chalopin_chepoi_hirai_osajda}
J{\'e}r{\'e}mie Chalopin, Victor Chepoi, Hiroshi Hirai, and Damian Osajda.
\newblock Weakly modular graphs and nonpositive curvature.
\newblock {\em Mem. Amer. Math. Soc.}, 2021.

\bibitem[CD95a]{charney_davis}
Ruth Charney and Michael~W. Davis.
\newblock {Finite {$K(\pi, 1)$}s for {A}rtin groups}.
\newblock In {\em Prospects in topology ({P}rinceton, {NJ}, 1994)}, volume 138
  of {\em Ann. of Math. Stud.}, pages 110--124. Princeton Univ. Press, 1995.

\bibitem[CD95b]{charney_davis_kpi1}
Ruth Charney and Michael~W. Davis.
\newblock The {$K(\pi,1)$}-problem for hyperplane complements associated to
  infinite reflection groups.
\newblock {\em J. Amer. Math. Soc.}, 8(3):597--627, 1995.

\bibitem[Cha]{charney_problems}
Ruth Charney.
\newblock {Problems related to Artin groups}.
\newblock {American Institute of Mathematics, http://people.brandeis.edu/\tild
  charney/papers/\verb|_|probs.pdf}.

\bibitem[CKS23]{ciobotaru_kramer_schwer_polyehdral_I}
Corina Ciobotaru, Linus Kramer, and Petra Schwer.
\newblock Polyhedral compactifications, {I}.
\newblock {\em Adv. Geom.}, 23(3):413--436, 2023.

\bibitem[CM09a]{caprace_monod_cat0_discrete}
Pierre-Emmanuel Caprace and Nicolas Monod.
\newblock Isometry groups of non-positively curved spaces: discrete subgroups.
\newblock {\em J. Topol.}, 2(4):701--746, 2009.

\bibitem[CM09b]{caprace_monod_CAT0_structure}
Pierre-Emmanuel Caprace and Nicolas Monod.
\newblock Isometry groups of non-positively curved spaces: structure theory.
\newblock {\em J. Topol.}, 2(4):661--700, 2009.

\bibitem[CMV23]{cumplido_martin_vaskou_parabolic}
Mar\'ia Cumplido, Alexandre Martin, and Nicolas Vaskou.
\newblock Parabolic subgroups of large-type {A}rtin groups.
\newblock {\em Math. Proc. Cambridge Philos. Soc.}, 174(2):393--414, 2023.

\bibitem[CN05]{chatterji_niblo}
Indira Chatterji and Graham Niblo.
\newblock From wall spaces to {$\rm CAT(0)$} cube complexes.
\newblock {\em Internat. J. Algebra Comput.}, 15(5-6):875--885, 2005.

\bibitem[CP05]{crisp_paris_representations_braid}
John Crisp and Luis Paris.
\newblock Representations of the braid group by automorphisms of groups,
  invariants of links, and {G}arside groups.
\newblock {\em Pacific J. Math.}, 221(1):1--27, 2005.

\bibitem[CP11]{corran_new_garside_braid}
Ruth Corran and Matthieu Picantin.
\newblock A new {G}arside structure for the braid groups of type {$(e,e,r)$}.
\newblock {\em J. Lond. Math. Soc. (2)}, 84(3):689--711, 2011.

\bibitem[CS11]{caprace_sageev}
Pierre-Emmanuel Caprace and Michah Sageev.
\newblock Rank rigidity for {CAT}(0) cube complexes.
\newblock {\em Geom. Funct. Anal.}, 21(4):851--891, 2011.

\bibitem[Dav15]{davis_coxeter}
Michael~W. Davis.
\newblock The geometry and topology of {C}oxeter groups.
\newblock In {\em Introduction to modern mathematics}, volume~33 of {\em Adv.
  Lect. Math. (ALM)}, pages 129--142. Int. Press, Somerville, MA, 2015.

\bibitem[Deh15]{garside}
Patrick Dehornoy.
\newblock {\em Foundations of {G}arside theory}, volume~22 of {\em EMS Tracts
  in Mathematics}.
\newblock European Mathematical Society (EMS), Z\"urich, 2015.
\newblock With Fran\c{c}ois Digne, Eddy Godelle, Daan Krammer and Jean Michel,
  Contributor name on title page: Daan Kramer.

\bibitem[Del72]{deligne}
Pierre Deligne.
\newblock Les immeubles des groupes de tresses g\'en\'eralis\'es.
\newblock {\em Invent. Math.}, 17:273--302, 1972.

\bibitem[Des16]{descombes_asymptotic_rank}
Dominic Descombes.
\newblock Asymptotic rank of spaces with bicombings.
\newblock {\em Math. Z.}, 284(3-4):947--960, 2016.

\bibitem[DL15]{descombes_lang_hyperbolicity}
Dominic Descombes and Urs Lang.
\newblock Convex geodesic bicombings and hyperbolicity.
\newblock {\em Geom. Dedicata}, 177:367--384, 2015.

\bibitem[DL16]{descombes_lang_flats}
Dominic Descombes and Urs Lang.
\newblock Flats in spaces with convex geodesic bicombings.
\newblock {\em Anal. Geom. Metr. Spaces}, 4(1):68--84, 2016.

\bibitem[DMW20]{dougherty_mccammond_witzel}
Michael Dougherty, Jon McCammond, and Stefan Witzel.
\newblock Boundary braids.
\newblock {\em Algebr. Geom. Topol.}, 20(7):3505--3560, 2020.

\bibitem[DPS22]{delucchi_paolini_salvetti_rank_3}
Emanuele Delucchi, Giovanni Paolini, and Mario Salvetti.
\newblock {Dual structures on Coxeter and Artin groups of rank three}.
\newblock {\em arXiv preprint arXiv:2206.14518}, 2022.

\bibitem[Duc18]{duchesne_CAT0_survey}
Bruno Duchesne.
\newblock Groups acting on spaces of non-positive curvature.
\newblock In {\em Handbook of group actions. {V}ol. {III}}, volume~40 of {\em
  Adv. Lect. Math. (ALM)}, pages 101--141. Int. Press, Somerville, MA, 2018.

\bibitem[FL05]{farrell_lafont_EZ}
F.~T. Farrell and J.-F. Lafont.
\newblock E{Z}-structures and topological applications.
\newblock {\em Comment. Math. Helv.}, 80(1):103--121, 2005.

\bibitem[FO20]{fukaya_oguni}
Tomohiro Fukaya and Shin-ichi Oguni.
\newblock A coarse {C}artan-{H}adamard theorem with application to the coarse
  {B}aum-{C}onnes conjecture.
\newblock {\em J. Topol. Anal.}, 12(3):857--895, 2020.

\bibitem[GJT98]{guivarch}
Yves Guivarc'h, Lizhen Ji, and J.~C. Taylor.
\newblock {\em {Compactifications of symmetric spaces\!}}
\newblock {Progr.~Math. {\bf 156}}. {Birkh\"auser}, 1998.

\bibitem[GP12]{godelle_paris}
Eddy Godelle and Luis Paris.
\newblock Basic questions on {A}rtin-{T}its groups.
\newblock In {\em Configuration spaces}, volume~14 of {\em CRM Series}, pages
  299--311. Ed. Norm., Pisa, 2012.

\bibitem[Gro87]{gromov_hyperbolic_groups}
M.~Gromov.
\newblock Hyperbolic groups.
\newblock In {\em Essays in group theory}, volume~8 of {\em Math. Sci. Res.
  Inst. Publ.}, pages 75--263. Springer, New York, 1987.

\bibitem[Hae21]{haettel_artin_cubic}
Thomas Haettel.
\newblock {Virtually cocompactly cubulated Artin-Tits groups}.
\newblock {\em Int. Math. Res. Not. IMRN}, (4):2919--2961, 2021.

\bibitem[Hae22a]{haettel_injective_buildings}
Thomas Haettel.
\newblock {Injective metrics on buildings and symmetric spaces}.
\newblock {\em Bull. Lond. Math. Soc.}, 2022.

\bibitem[Hae22b]{haettel_bourbaki_paolini_salvetti}
Thomas Haettel.
\newblock {La conjecture du $K(\pi,1)$ pour les groupes d'Artin affines
  (d'apr{\`e}s Paolini et Salvetti)}.
\newblock {\em {S\'em. Bourbaki. Exp. No. 1195. Ast\'erisque}}, 438, 2022.

\bibitem[Hae24]{haettel_helly_kpi1}
Thomas Haettel.
\newblock {Lattices, injective metrics and the $K(\pi,1)$ conjecture}.
\newblock {\em Algebr. Geom. Topol.}, 2024.

\bibitem[HH23]{haettel_huang_weakly_modular}
Thomas Haettel and Jingyin Huang.
\newblock {Lattices, Garside structures and weakly modular graphs}.
\newblock {\em J. Algebra}, 2023.

\bibitem[HHP21]{haettel_hoda_petyt}
Thomas Haettel, Nima Hoda, and Harry Petyt.
\newblock {Coarse injectivity, hierarchical hyperbolicity, and
  semihyperbolicity}.
\newblock {\em to appear in Geom. Topol.}, 2021.

\bibitem[HHP23]{haettel_hoda_petyt_Lp_metrics}
Thomas Haettel, Nima Hoda, and Harry Petyt.
\newblock $\ell^p$ metrics on cell complexes.
\newblock {\em arXiv preprint arXiv:2302.03801}, 2023.

\bibitem[Hir20]{hirai_uniform_modular}
Hiroshi Hirai.
\newblock Uniform modular lattices and affine buildings.
\newblock {\em Adv. Geom.}, 20(3):375--390, 2020.

\bibitem[Hir21]{hirai_CAT0}
Hiroshi Hirai.
\newblock A nonpositive curvature property of modular semilattices.
\newblock {\em Geom. Dedicata}, 214:427--463, 2021.

\bibitem[HKS16]{b6}
Thomas Haettel, Dawid Kielak, and Petra Schwer.
\newblock The 6-strand braid group is {${\rm CAT}(0)$}.
\newblock {\em Geom. Dedicata}, 182:263--286, 2016.

\bibitem[HM12]{hatcher_margalit_torelli}
Allen Hatcher and Dan Margalit.
\newblock Generating the {T}orelli group.
\newblock {\em Enseign. Math. (2)}, 58(1-2):165--188, 2012.

\bibitem[HO21]{huang_osajda_helly}
Jingyin Huang and Damian Osajda.
\newblock {Helly meets Garside and Artin}.
\newblock {\em Invent. Math.}, 225(2):395--426, 2021.

\bibitem[HSWW17]{haettel_schilling_wienhard}
Thomas Haettel, Anna-Sofie Schilling, Cormac Walsh, and Anna Wienhard.
\newblock {Horofunction compactifications of symmetric spaces}.
\newblock 2017.
\newblock {arXiv:1705:05026}.

\bibitem[HW08]{haglund_wise_ccc}
Fr\'ed\'eric Haglund and Daniel~T. Wise.
\newblock Special cube complexes.
\newblock {\em Geom. Funct. Anal.}, 17(5):1551--1620, 2008.

\bibitem[HW12]{haglund_wise}
Fr{\'e}d{\'e}ric Haglund and Daniel~T. Wise.
\newblock A combination theorem for special cube complexes.
\newblock {\em Ann. of Math. (2)}, 176(3):1427--1482, 2012.

\bibitem[Isb64]{isbell}
J.~R. Isbell.
\newblock Six theorems about injective metric spaces.
\newblock {\em Comment. Math. Helv.}, 39:65--76, 1964.

\bibitem[Jeo23]{jeong}
Seong~Gu Jeong.
\newblock The seven-strand braid group is {${\rm CAT}(0)$}.
\newblock {\em Manuscripta Math.}, 171(3-4):563--581, 2023.

\bibitem[JS06]{januskiewicz_swiatkowski_systolic}
Tadeusz Januszkiewicz and Jacek \'{S}wi\c{a}tkowski.
\newblock Simplicial nonpositive curvature.
\newblock {\em Publ. Math. Inst. Hautes \'{E}tudes Sci.}, (104):1--85, 2006.

\bibitem[JS17]{ji_schilling_horofunction_polyhedral}
Lizhen Ji and Anna-Sofie Schilling.
\newblock Toric varieties vs. horofunction compactifications of polyhedral
  norms.
\newblock {\em Enseign. Math.}, 63(3-4):375--401, 2017.

\bibitem[Kak92]{kakimizu}
Osamu Kakimizu.
\newblock Incompressible spanning surfaces and maximal fibred submanifolds for
  a knot.
\newblock {\em Math. Z.}, 210(2):207--223, 1992.

\bibitem[KL20]{kleiner_lang_higher_rank_hyperbolicity}
Bruce Kleiner and Urs Lang.
\newblock Higher rank hyperbolicity.
\newblock {\em Invent. Math.}, 221(2):597--664, 2020.

\bibitem[KPP12]{Kpoczynski_pak_przytycki_acute_triangulations}
Eryk Kopczy\'{n}ski, Igor Pak, and Piotr Przytycki.
\newblock Acute triangulations of polyhedra and {${\Bbb R}^N$}.
\newblock {\em Combinatorica}, 32(1):85--110, 2012.

\bibitem[KR17]{kasprowski_rueping}
Daniel Kasprowski and Henrik R\"uping.
\newblock The {F}arrell-{J}ones conjecture for hyperbolic and {CAT}(0)-groups
  revisited.
\newblock {\em J. Topol. Anal.}, 9(4):551--569, 2017.

\bibitem[Lan13]{lang}
Urs Lang.
\newblock Injective hulls of certain discrete metric spaces and groups.
\newblock {\em J. Topol. Anal.}, 5(3):297--331, 2013.

\bibitem[Lea13]{leary_CCC}
Ian~J. Leary.
\newblock A metric {K}an-{T}hurston theorem.
\newblock {\em J. Topol.}, 6(1):251--284, 2013.

\bibitem[LW11]{lemmens_walsh_isometries_hilbert}
Bas Lemmens and Cormac Walsh.
\newblock Isometries of polyhedral {H}ilbert geometries.
\newblock {\em J. Topol. Anal.}, 3(2):213--241, 2011.

\bibitem[LY21]{leuzinger_young}
Enrico Leuzinger and Robert Young.
\newblock Filling functions of arithmetic groups.
\newblock {\em Ann. of Math. (2)}, 193(3):733--792, 2021.

\bibitem[McC09]{mccammond_CAT0_survey}
Jon McCammond.
\newblock Constructing non-positively curved spaces.
\newblock In {\em Geometric and cohomological methods in group theory}, volume
  358 of {\em London Math. Soc. Lecture Note Ser.}, pages 162--224. Cambridge
  Univ. Press, Cambridge, 2009.

\bibitem[McC15]{mccammond_failure_lattice_property}
Jon McCammond.
\newblock Dual euclidean {A}rtin groups and the failure of the lattice
  property.
\newblock {\em J. Algebra}, 437:308--343, 2015.

\bibitem[McC17]{mccammond_mysterious}
Jon McCammond.
\newblock {The mysterious geometry of Artin groups}.
\newblock {2017}.
\newblock {http://web.math.ucsb.edu/\tild
  jon.mccammond/papers/mysterious-geometry.pdf}.

\bibitem[Mie14]{miesch_CCC}
Benjamin Miesch.
\newblock Injective metrics on cube complexes.
\newblock {\em arXiv preprint arXiv:1411.7234}, 2014.

\bibitem[Mie17]{miesch}
Benjamin Miesch.
\newblock The {C}artan-{H}adamard theorem for metric spaces with local geodesic
  bicombings.
\newblock {\em Enseign. Math.}, 63(1-2):233--247, 2017.

\bibitem[MM99]{masur_minsky}
Howard~A. Masur and Yair~N. Minsky.
\newblock Geometry of the complex of curves. {I}. {H}yperbolicity.
\newblock {\em Invent. Math.}, 138(1):103--149, 1999.

\bibitem[MS17]{mccammond_sulway}
Jon McCammond and Robert Sulway.
\newblock Artin groups of {E}uclidean type.
\newblock {\em Invent. Math.}, 210(1):231--282, 2017.

\bibitem[NR98]{niblo_roller}
Graham~A. Niblo and Martin~A. Roller.
\newblock Groups acting on cubes and {K}azhdan's property ({T}).
\newblock {\em Proc. Amer. Math. Soc.}, 126(3):693--699, 1998.

\bibitem[NR03]{niblo_reeves}
G.~A. Niblo and L.~D. Reeves.
\newblock Coxeter groups act on {${\rm CAT}(0)$} cube complexes.
\newblock {\em J. Group Theory}, 6(3):399--413, 2003.

\bibitem[Nus88]{nussbaum_hilbert_metric}
Roger~D. Nussbaum.
\newblock Hilbert's projective metric and iterated nonlinear maps.
\newblock {\em Mem. Amer. Math. Soc.}, 75(391):iv+137, 1988.

\bibitem[OP09]{osajda_przytycki_systolic_boundary}
Damian Osajda and Piotr Przytycki.
\newblock Boundaries of systolic groups.
\newblock {\em Geom. Topol.}, 13(5):2807--2880, 2009.

\bibitem[OV24]{osajda_valiunas}
Damian Osajda and Motiejus Valiunas.
\newblock Helly groups, coarsely {H}elly groups, and relative hyperbolicity.
\newblock {\em Trans. Amer. Math. Soc.}, 377(3):1505--1542, 2024.

\bibitem[Pao21]{paolini_dual_approach}
Giovanni Paolini.
\newblock {The dual approach to the $K(\pi,1)$ conjecture}.
\newblock {\em arXiv preprint arXiv:2112.05255}, 2021.

\bibitem[Par14]{paris_kpi1}
Luis Paris.
\newblock {$K(\pi,1)$} conjecture for {A}rtin groups.
\newblock {\em Ann. Fac. Sci. Toulouse Math. (6)}, 23(2):361--415, 2014.

\bibitem[Pic22]{picantin_amalgams_Garside}
Matthieu Picantin.
\newblock Cyclic amalgams, {HNN} extensions, and {G}arside one-relator groups.
\newblock {\em J. Algebra}, 607(part B):437--465, 2022.

\bibitem[Prz08]{przytycki_systolic_fixed_point}
Piotr Przytycki.
\newblock The fixed point theorem for simplicial nonpositive curvature.
\newblock {\em Math. Proc. Cambridge Philos. Soc.}, 144(3):683--695, 2008.

\bibitem[PS12]{przytycki_schultens_kakimizu}
Piotr Przytycki and Jennifer Schultens.
\newblock Contractibility of the {K}akimizu complex and symmetric {S}eifert
  surfaces.
\newblock {\em Trans. Amer. Math. Soc.}, 364(3):1489--1508, 2012.

\bibitem[PS21]{paolini_salvetti}
Giovanni Paolini and Mario Salvetti.
\newblock Proof of the {$K(\pi,1)$} conjecture for affine {A}rtin groups.
\newblock {\em Invent. Math.}, 224(2):487--572, 2021.

\bibitem[Rol98]{roller}
Martin Roller.
\newblock {Poc Sets, Median Algebras and Group Actions}.
\newblock {arXiv:1607.07747}, 1998.

\bibitem[Ron89]{ronan}
M.A. Ronan.
\newblock {\em {Lectures on buildings\!}}
\newblock Persp. Math. {\bf 7}. Academic Press, 1989.

\bibitem[Sag95]{sageev}
Michah Sageev.
\newblock Ends of group pairs and non-positively curved cube complexes.
\newblock {\em Proc. London Math. Soc. (3)}, 71(3):585--617, 1995.

\bibitem[Sch10]{schultens_kakimizu}
Jennifer Schultens.
\newblock The {K}akimizu complex is simply connected.
\newblock {\em J. Topol.}, 3(4):883--900, 2010.
\newblock With an appendix by Michael Kapovich.

\bibitem[ST88]{scharlemann_thompson_kakimizu}
Martin Scharlemann and Abigail Thompson.
\newblock Finding disjoint {S}eifert surfaces.
\newblock {\em Bull. London Math. Soc.}, 20(1):61--64, 1988.

\bibitem[Wah13]{wahl}
Nathalie Wahl.
\newblock Homological stability for mapping class groups of surfaces.
\newblock In {\em Handbook of moduli. {V}ol. {III}}, volume~26 of {\em Adv.
  Lect. Math. (ALM)}, pages 547--583. Int. Press, Somerville, MA, 2013.

\bibitem[Wal07]{walsh_horofunction_normed_space}
Cormac Walsh.
\newblock The horofunction boundary of finite-dimensional normed spaces.
\newblock {\em Math. Proc. Cambridge Philos. Soc.}, 142(3):497--507, 2007.

\bibitem[Wal08]{walsh_horofunction_hilbert}
Cormac Walsh.
\newblock The horofunction boundary of the {H}ilbert geometry.
\newblock {\em Adv. Geom.}, 8(4):503--529, 2008.

\bibitem[Web20]{webb_non_CAT0}
Richard C.~H. Webb.
\newblock Contractible, hyperbolic but non-{${\rm CAT}(0)$} complexes.
\newblock {\em Geom. Funct. Anal.}, 30(5):1439--1463, 2020.

\bibitem[Wen05]{wenger}
S.~Wenger.
\newblock Isoperimetric inequalities of {E}uclidean type in metric spaces.
\newblock {\em Geom. Funct. Anal.}, 15(2):534--554, 2005.

\bibitem[Wis21]{wise_book_hierarchy}
Daniel~T. Wise.
\newblock {\em The structure of groups with a quasiconvex hierarchy}, volume
  209 of {\em Annals of Mathematics Studies}.
\newblock Princeton University Press, Princeton, NJ, [2021] \copyright 2021.

\end{thebibliography}

\end{document}